\documentclass[a4paper,12pt]{article}

\usepackage[english,francais]{babel}
\usepackage[T1]{fontenc}
\usepackage{lmodern}
\usepackage[hmargin=2.5cm,vmargin=30mm]{geometry}
\usepackage{amsthm,amssymb,amsmath}
\usepackage{url}
\usepackage{wrapfig}

\def\NN{\mathbb{N}}
\def\ZZ{\mathbb{Z}}
\def\QQ{\mathbb{Q}}
\def\RR{\mathbb{R}}
\def\CC{\mathbb{C}}
\def\PP{\mathbb{P}}

\def\HHH{\mathcal{H}}

\def\LLL{\mathcal{L}}
\def\RRR{\mathcal{R}}

\def\Box{\operatorname{Box}}
\def\SL{\operatorname{SL}}

\def\Aff{\operatorname{Aff}}

\def\Deck{\operatorname{Deck}}
\def\Hdim{\operatorname{Hdim}}

\renewcommand{\Re}{\operatorname{Re}}
\renewcommand{\Im}{\operatorname{Im}}

\def\Xsurf{\mathrm{X}}
\def\Lsurf{\mathrm{L}}

\newtheorem{proposition}{Proposition}
\newtheorem{theorem}{Theorem}
\newtheorem{corollary}{Corollary}
\newtheorem{lemma}{Lemma}
\newtheorem{definition}{Definition}

\usepackage{graphicx}
\def\picinput#1{\includegraphics{#1.pdf}}

\author{V. Delecroix}
\title{Divergent trajectories in the periodic wind-tree model}
\date{2011}

\begin{document}

\maketitle

\selectlanguage{english}

\begin{abstract}
The periodic wind-tree model is a family $T(a,b)$ of billiards in the plane in which identical rectangular scatterers of size $a \times b$ are disposed at each integer point. It was proven by P.~Hubert, S.~Leli\`evre and S.~Troubetzkoy that for a residual set of parameters $(a,b)$ the billiard flow in $T(a,b)$ is recurrent in almost every direction. We prove that for many parameters $(a,b)$ there exists a set $\Lambda \subset S^1$ of positive Hausdorff dimension such that for every $\theta \in \Lambda$ every billiard trajectory in $T(a,b)$ with initial angle $\theta$ is divergent.

\vskip 0.5\baselineskip

\selectlanguage{francais}
\begin{center}\textbf{R\'esum\'e}\end{center}
\vskip 0.5\baselineskip
{\bf Trajectoires divergentes pour vent dans les arbres}

\vspace{0.5\baselineskip}
Le ``vent dans les arbres'' est une famille de billards infinis $T(a,b)$ d\'efinis de la mani\`ere suivante. Dans le plan euclidien $\RR^2$, on place des rectangles de taille $a \times b$ \`a chaque point entier. Une particule (identifi\'ee \`a un point) se d\'eplace en ligne droite et rebondit de mani\`ere \'elastique sur les obstacles. P.~Hubert, S.~Leli\`evre et S.~Troubetzkoy ont d\'emontr\'e qu'il existait un $G_\delta$ dense de param\`etres $(a,b)$ pour lesquels, dans presque toute direction $\theta \in S^1$, le flot du billard $T(a,b)$ dans la direction $\theta$ est r\'ecurrent. Nous prouvons que pour certains param\`etres $(a,b)$, il existe un ensemble $\Lambda \subset S^1$ de mesure de Hausdorff positive tel que pour tout $\theta \in \Lambda$ toute trajectoire dans le billard $T(a,b)$ dont l'angle de d\'epart est $\theta$ est divergente.
\end{abstract}

\selectlanguage{english}

\section{Introduction}
\label{section:introduction}
\begin{figure}[!ht]
\begin{center}
\includegraphics[scale=0.75]{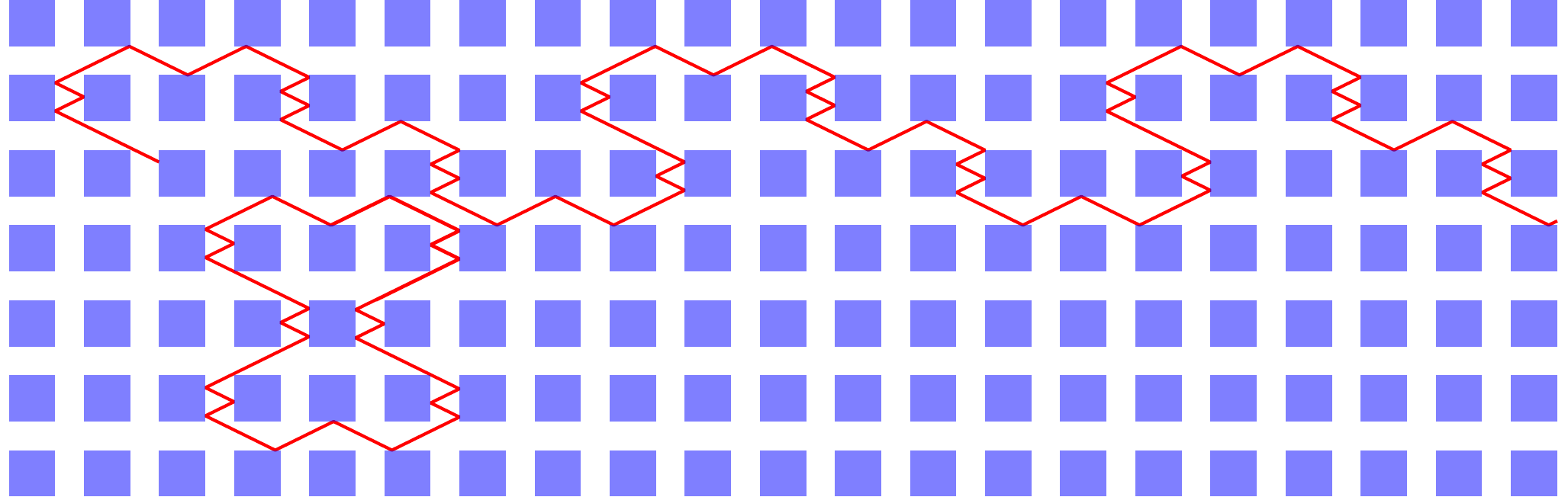}
\end{center}
\caption{A periodic and a divergent orbits in a half-divergent direction of the ``golden wind-tree model'' with parameters $a=b=\Phi=\frac{\sqrt{5}-1}{2}$.}
\label{fig:golden_examples}
\end{figure}

We study periodic versions of the wind-tree model introduced by P.~Ehrenfest and T.~Ehrenfest in 1912~\cite{EhrenfestEhrenfest12}. A point (``wind'') moves in the plane and collides with rectangular scatterers (``trees'') with the usual law of reflexion. In the periodic version of the wind-tree model, due to J.~Hardy and J.~Weber~\cite{HardyWeber80}, the scatterers are identical rectangular obstacles located periodically in the plane, every obstacle centered at each point of $\ZZ^2$. The scatterers are rectangles of size $a \times b$, with $0 < a < 1$, $0 < b < 1$. We denote by $T(a,b)$ the subset of the plane obtained by removing the obstacles and name its billiard \emph{the wind-tree model}. Our aim is to understand some of its dynamical properties (see Figure~\ref{fig:golden_examples} for two different behaviors in the golden wind-tree table $T(\Phi,\Phi)$).

The phase space of the billiard is naturally $T(a,b) \times S^1$. Each barrier in $T(a,b)$ is either horizontal or vertical. Hence, for the point $(x,\theta) \in T(a,b) \times S^1$ and for every time $t$, the possible angles for the orbit of $(x,\theta)$ in $T(a,b)$ at time $t$ are $\theta$, $-\theta$, $\pi/2-\theta$  or $-\pi/2+\theta$. Let $\tau_h: \theta \mapsto -\theta$ and $\tau_v: \theta \mapsto \pi/2-\theta$ be respectively the horizontal and vertical reflexions and $K = \langle \tau_h, \tau_v \rangle \simeq \ZZ/2 \times \ZZ/2$ the group they generate. We define the \emph{billiard flow in direction $\theta$} in $T(a,b)$ to be the map $\phi^\theta_t: T(a,b) \times K \rightarrow T(a,b) \times K$ which is defined as follows. Let $(x_0,\kappa_0)$ be an element of $T(a,b) \times K$, then $(x_t,\kappa_t) = \phi_t^\theta(x,\kappa)$ is such that if a ball has an initial position $x_0$ and an initial angle $\kappa_0(\theta)$ then after time $t$ it has position $x_t$ and direction $\kappa_t(\theta)$. We will often consider the quantity $\phi_t^\theta(x,\tau)$ as an element of $T(a,b)$ and write $\phi_t^\theta(x)$ for $\phi_t^\theta(x,id)$.

Let $d$ be the Euclidean distance in $\RR^2$. We say that the flow in direction $\theta$ is \emph{recurrent}, if for almost all points $x$ in $T(a,b)$ we have $\displaystyle \liminf_{t \to \infty} d(x,\phi^\theta_t(x)) = 0$. We say that it is {\em divergent} if for almost all points $x \in T(a,b)$ $\displaystyle \liminf_{t \to \infty} d(x,\phi_t^\theta(x)) = +\infty$. P.~Hubert, S.~Leli\`evre and S.~Troubetkoy \cite{HubertLelievreTroubetzkoy} exhibit a residual set $\mathcal{E} \subset (0,1) \times (0,1)$ such that for any parameters $(a,b) \in \mathcal{E}$ for almost all $\theta \in S^1$ (with respect to the Lebesgue measure on $S^1$) the flow in $T(a,b)$ in direction $\theta$ is recurrent. In the present paper we study the opposite behavior: the set of parameters $(a,b,\theta)$ for which the flow in $T(a,b)$ in direction $\theta$ is divergent. As a consequence of our main result (Theorem~\ref{thm:even_implies_div}) we obtain
\begin{theorem} \label{thm:divergent_parameters}
If $a$ and $b$ are either rational or quadratic of the form $1/(1-a) = x + y \sqrt{D}$ and $1/(1-b) = (1-x)+y \sqrt{D}$ with $x$ and $y$ rationals and $D > 1$ a square-free integer, then there exists a dense set $\Lambda \subset S^1$ of Hausdorff dimension not smaller than $1/2$ such that for every $\theta \in \Lambda$ and every point $x$ in $T(a,b)$ with infinite forward orbit $\liminf d(x,\phi_t^\theta(x)) = \infty$. In particular the flow $\phi_t^\theta$ is divergent.
\end{theorem}
The subset $\Lambda$ which appears in the above theorem is made explicit in Proposition~\ref{prop:infinitely_renormalizable_veech}.

The strategy to prove Theorem~\ref{thm:divergent_parameters} is very similar to the one used in \cite{HubertLelievreTroubetzkoy}. We explain the idea for the special case $T(1/2,1/2)$. For every angle $\theta \in (0,\pi/2)$ for which the slope $\tan \theta$ is rational, the billiard flow in $T(a,b)$ in direction $\theta$ has a periodic behavior. Two important types of slopes are of interest for our purpose: \emph{half-divergent slope} and \emph{periodic slope} (see Figure~\ref{fig:golden_examples} for an example of a half-divergent slope). We explain the definition on two examples. In the horizontal direction $\theta=0$ in $T(1/2,1/2)$, there is a bunch of trajectories which reflect between two consecutive scatterers spaced by $(1,0)$ while the others go to infinity: $0$ is a half-divergent slope. On the contrary, in the direction $\theta = \pi/4$, all trajectories are periodic with the same period: the slope $\tan(\pi/4) = 1$ is periodic. To prove recurrence, the strategy of \cite{HubertLelievreTroubetzkoy} consists in approximate a generic slopes by rational ones which correspond to directions of periodic type in $T(1/2,1/2)$. To build divergent trajectories in the same billiard table, we use slopes are in a sense badly approximate by slopes of periodic type.

The proof of our main result uses a renormalization algorithm due to S.~Ferenczi and L.~Zamboni~\cite{FerencziZamboni10_struct, FerencziZamboni10_eig}. Their induction operate on interval exchange transformations and we give a geometric interpretation on translation surfaces using suspensions. Similar geometric renormalization is described by C.~Ulcigra\"i and J.~Smillie for the regular ``octagon''~\cite{SmillieUlcigrai}. The geometric interpretation we use was known in greater generality by P.~Hubert and C.~Ulcigra\"i~\cite{HubertUlcigrai}.

We first consider a discretization of the flow $\phi^\theta_t$ in $T(a,b)$ and prove that the distance $d(x,\phi_t^\theta(x))$ corresponds to a Birkhoff sum of a function over an interval exchange transformation $g = g_{a,b,\theta}$. Then we build a set of parameters $(a,b,\theta)$ by imposing some conditions in the Ferenczi-Zamboni induction of $g$. For those parameters, we have a very simple continued fraction algorithm-like which is define as follows. For a 4-tuple of positive real numbers $Z=(x_1,x_2,y_1,y_2)$ define
\[
F(Z) = \left(
y_1,\ y_2,\ 
x_1 - \left\lfloor \frac{x_1}{y_1+y_2} \right\rfloor (y_1 + y_2),\ 
x_2 - \left\lfloor \frac{x_2}{y_2} \right\rfloor y_2 
\right)
\]
and 
\[
m(Z) = \left\lfloor \frac{x_1}{y_1+y_2} \right\rfloor
\quad \text{and} \quad
n(Z) = \left\lfloor \frac{x_2}{y_2} \right\rfloor
\]
where $\lfloor . \rfloor$ designs the floor. If $Z$ satisfies
\begin{equation} \label{eq:non_renormalizable}
x_1 + x_2 > y_1 > x_2 \quad \text{and} \quad y_1 + y_2 > x_1 > y_2
\end{equation} 
then $F^2(Z) = Z$. We say that the quadruple $Z$ is $F$-\emph{renormalizable} if for all $k \geq 0$, $F^k(Z)$ does not satisfy~(\ref{eq:non_renormalizable}). To a $F$-renormalizable quadruple $Z$ we associate an infinite sequence of $2$-tuples $((m_k,n_k))_{k \in \NN}$ defined by $m_k(Z) = m(F^k(Z))$ and $n_k(Z) = n(F^k(Z))$. We call the sequence $((m_k(Z),n_k(Z))_{k \geq 0}$ the $F$-\emph{convergents} of $Z$. The set of $F$-renormalizable quadruples defines an uncountable set of zero Lebesgue measure.
\begin{proposition}
If $Z=(x_1,x_2,y_1,y_2)$ is $F$-renormalizable then
\begin{itemize}
\item for $k \geq 0$, if $m_k = 0$ then $m_{k+1} \not= 0$ and $n_{k+1} = 0$,
\item for infinitely many $k$, $m_{2k} \not= 0$. The same is true for $m_{2k+1}$, $n_{2k}$ and $n_{2k+1}$.
\end{itemize}
Conversely, if $((m_k,n_k))_{k \geq 0}$  is a sequence of 2-tuples of non negative integers that satisfy the above condition, then there exists at least one quadruple $Z$ such that for all $k$, $m_k(Z) = m_k$ and $n_k(Z) = n_k$.
\end{proposition}
Using the above description of $F$-renormalizable slopes, our main result is
\begin{theorem} \label{thm:even_implies_div}
Let $Z = (x_1,x_2,y_1,y_2) \in \RR^4_+$ and $(a,b,\theta) \in (0,1)^2 \times [0,\pi/2)$ be related by
\[
x_1 = (1-b) \cos \theta \quad  x_2 = b \cos \theta \quad y_1 = (1-a) \sin \theta \quad  \text{and} \quad y_2 = a \sin \theta.
\]
Assume that $Z$ is $F$-renormalizable and let $((m_k,n_k))_{k \geq 0}$ be its $F$-convergents. If for all $k$, $n_k \equiv 0 \mod 2$, then any infinite forward trajectory in direction $\theta$ in $T(a,b)$ is self-avoiding. In particular the flow in direction $\theta$ in $T(a,b)$ is divergent.
\end{theorem}

\vspace{0.5cm}

The infinite billiard $T(a,b)$ can be considered as a particular case of $\ZZ^2$-periodic translation surface (with finite quotient). For $\ZZ$-periodic translation surfaces the recurrence of the flow follows from general results on $1$-dimensional cocycles. For highly symmetric examples, the ergodicity of $\ZZ$-periodic translation surface the flow can be established (see \cite{HubertSchmithuesen10}, \cite{HubertWeiss}, \cite{HooperHubertWeiss}). The main difficulty of the wind-tree model comes from dimension $2$.

We mention other results on the wind-tree-model. As we explained above, it is proven in~\cite{HubertLelievreTroubetzkoy} that for a residual set of parameters $(a,b)$ for almost all angles $\theta$ the flow $\phi^\theta_t$ in $T(a,b)$ is recurrent. The problem of diffusion is studied in~\cite{DelecroixHubertLelievre11} (it is proven that for almost all parameters and almost all angles $\theta$ the polynomial growth of $d(p,\phi^\theta_t(p))$ is $2/3$). The ergodic decomposition for irrational parameters $(a,b)$ and some angles $\theta$ with rational slopes is done in~\cite{ConzeGutkin}. We would also mention that the main result in~\cite{Hooper} gives the following result
\begin{theorem}[\cite{Hooper}]
If $a=p/q$ and $b=r/s$ are rationals with $p$, $r$ odd and $r$, $s$ even, then there exists a dense set $\Lambda \subset S^1$ of Hausdorff dimension not less than $1/2$ such that for every $\theta \in S^1$ the billiard flow $\phi^\theta_t$ is ergodic.
\end{theorem}

\vspace{0.5cm}

The paper is organized as follows. In Section~\ref{section:translation_surfaces} we define translation surfaces and interval exchange transformations. We build the discretization of the flow $\phi^\theta_t$ as a $\ZZ^2$-cocycle over an interval exchange transformation. Next in Section~\ref{section:divergent_cocycle} we recall the Ferenczi-Zamboni induction and see the relation with the map $F$ defined above. The proof that Theorem~\ref{thm:even_implies_div} implies Theorem~\ref{thm:divergent_parameters} is left to Section~\ref{subsection:veech_parameters}. The proof uses the classification of Veech surfaces in genus $2$ by K.~Calta~\cite{Calta04} and C.~McMullen~\cite{McMullen03}. The proof of Theorem~\ref{thm:even_implies_div} is postponed to Section~\ref{subsection:main_proof}.

\vspace{0.5cm}

\textbf{Aknowledgments:} The author would like to thank Pascal Hubert and Samuel Leli\`evre for introducing him to the known results about the wind-tree model $T(p/q,r/s)$ and more generally to the theory of finite and infinite square tiled surfaces. Many experimentations have been done with the math software Sage \cite{Sage,Sage-Combinat}. The script (for computations and drawings) as well as a collection of pictures are available on the web page of the author.


\section{The wind-tree cocycle} \label{section:translation_surfaces}
In this section we build a discretization of the billiard flow in $T(a,b)$.

\subsection{Translation surface and Poincare maps of the linear flow} \label{subsection:translation_surfaces}
A \emph{flat surface} is a compact oriented surface $X$ endowed with a flat metric defined on $X \backslash \Sigma$ where $\Sigma \subset X$ is a finite set of points which are conic type singularities for the metric. It is a \emph{translation surface} if moreover the holonomy given by parallel transport in $X \backslash \Sigma$ is trivial. Concretely, any translation surface can be built from a finite set of polygons $P_i$ in $\RR^2$ and identifying pairs of edges with translations. We refer to the survey of A.~Zorich~\cite{Zorich06} and the notes of M.~Viana~\cite{Viana} for the latter construction and other equivalent definitions of translation surfaces.

Let $X$ be a translation surface, the absence of holonomy implies that directions are globally defined. The geodesic flow on the tangent bundle of $X$ preserves directions and can be defined on $X$ as soon as we specify the direction and the speed. We assume that in a translation surface a fixed direction is given which we call \emph{vertical}. The \emph{linear flow} of $X$ is the unit speed geodesic flow in the vertical direction on $X$. The flow in a rational billiards can be transformed into the linear flow of a translation surface. We will use this construction to define the wind-tree cocycle (see Section~\ref{subsection:windtree_cocycle}).

\begin{wrapfigure}{r}{6cm}
\begin{center}\picinput{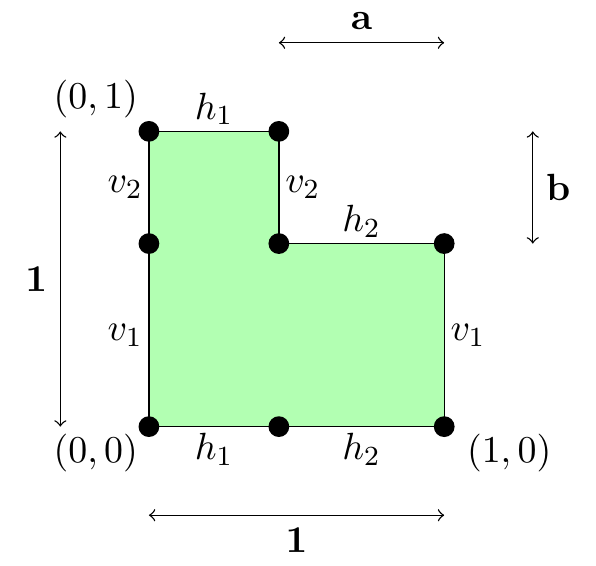}\end{center}
\caption{The surface $\Lsurf(a,b)$ (letters $h_i$ and $v_i$ indicate gluings).}
\label{fig:L_ab}
\end{wrapfigure}
In this paper we mainly focus on the family of translation surfaces $\Lsurf(a,b)$ where $0 <a < 1$ and $0 < b <1$ are two parameters. The surface $\Lsurf(a,b)$ is built as follows (see also Figure~\ref{fig:L_ab}). Take the polygon with vertices $(0,0)$, $(1,0)$, $(1,1-b)$, $(1-a,1-b)$, $(1-a,1)$, and $(0,1)$ and identify the following pairs of edges
\begin{itemize}
\item $[(0,0),(1-a,0)]$ with $[(0,1),(1-a,1)]$ (labelled $h_1$ in Figure~\ref{fig:L_ab}),
\item $[(0,0),(0,1-b)]$ with $[(1,0),(1,1-b)]$ (labelled $v_1$),
\item $[(1-a,0),(1,0)]$ with $[(1-a,1-b),(1,1-b)]$ (labelled $h_2$) and
\item $[(0,1-b),(0,1)]$ with $[(1-a,1-b),(1-a,1)]$ (labelled $v_2$).
\end{itemize}

The translation surface $\Lsurf(a,b)$ is a genus $2$ surface that has one conic singularity of angle $6\pi$. The equivalence classes of translation surfaces with a single conic singularity of angle $6\pi$ form a moduli space denoted $\mathcal{H}(2)$. The number $2$ in $\mathcal{H}(2)$ does not refer to the genus but to the \emph{degree} of the singularity. Any surface in $\mathcal{H}(2)$ can be built from a polygon which is L-shaped but non necessarily right-angled (see \cite{Calta04}, \cite{HubertLelievre06} and \cite{McMullen03}).

Let $X$ be a translation surface and $I \subset X$ an horizontal segment (or any segment transverse to the linear flow of $X$). The first return map on $I$ is defined for every point in $I$ for which the orbit under the linear flow returns to $I$ before reaching a singularity of the metric. There is a natural Lebesgue measure on $I$ induced from the flat metric of $X$ which is preserved by the first return map. Out of the discontinuities, the first return map on $I$ is a translation. For our purpose, it will be easier to work with Poincare maps obtained on transversal made of more than one interval.

\subsection{Interval exchange transformations and quadrangulations} \label{section:train_tracks_and_quadrangulations}
We mainly follows the presentation of \cite{FerencziZamboni10_eig,FerencziZamboni10_struct} and give a geometric point of view on their induction procedure. We use here the letters $\ell$ (for left) and $r$ (for right) whereas in \cite{FerencziZamboni10_eig,FerencziZamboni10_struct} the letters $m$ (for minus) and $p$ (for plus) are used.

\begin{figure}[!ht]
\begin{center}
\picinput{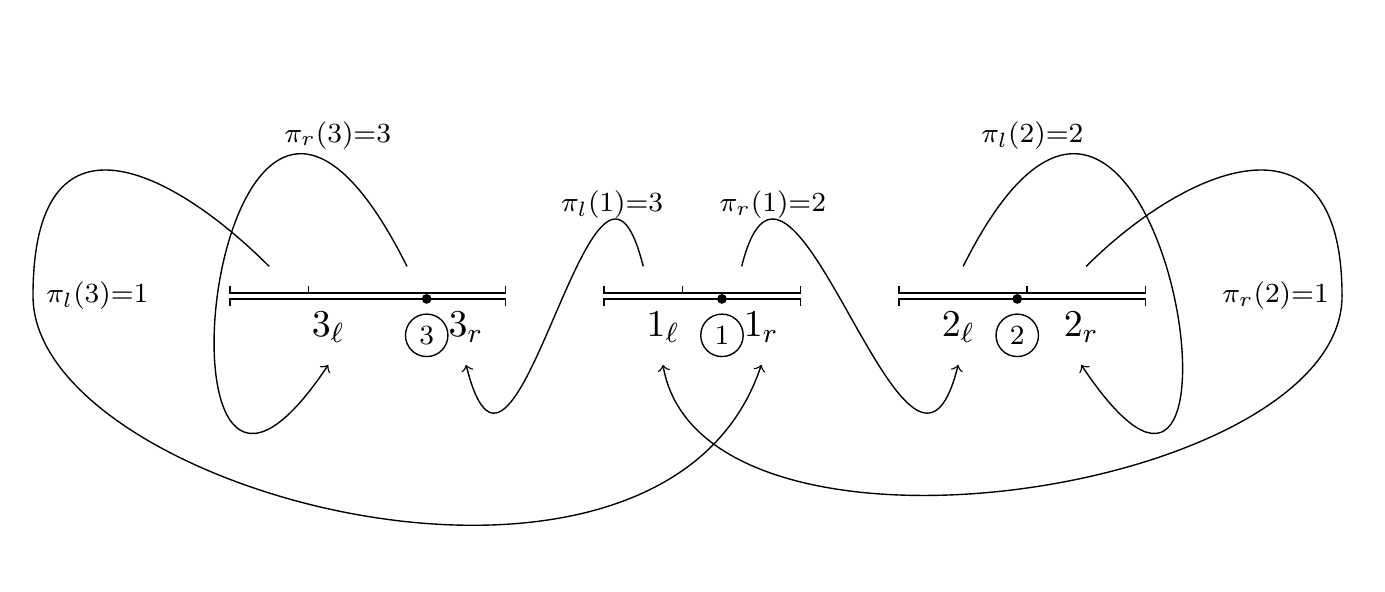}

\picinput{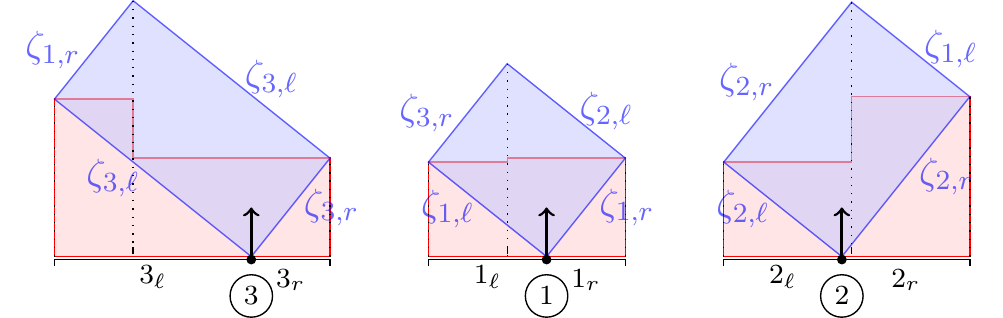}

\picinput{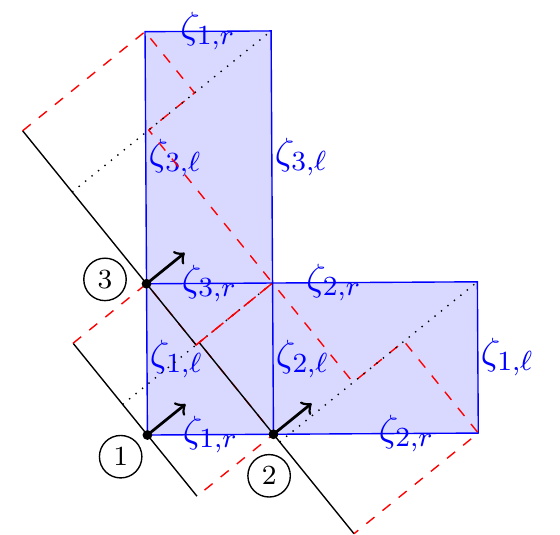}
\end{center}
\caption{An interval exchange transformation with permutation datum $\pi_\ell = (1,3)(2)$, $\pi_r = (1,2)(3)$ and two views of the unique zippered rectangles associated to it which gives a translation surface of the form $r_\theta \cdot \Lsurf(a,b)$.}
\label{fig:from_train_track_to_L}
\end{figure}

Let $\lambda = ((\lambda_{1,\ell},\lambda_{1,r}),\ldots,(\lambda_{d,\ell},\lambda_{d,r}))$ be a vector of $d$ pairs of real numbers where $\lambda_{i,\ell} < 0$ and $\lambda_{i,r} > 0$. Set $E_i = ]\lambda_{i,\ell},\lambda_{i,r}[$ for $i=1,\dots,d$. We define a map $T$ on the disjoint union $E = E_1 \sqcup E_2 \sqcup \dots \sqcup E_d$ for which the origin of every interval $E_i$ is a discontinuity of $T^{-1}$. The combinatorics of $T$ depends on a pair of permutations $\pi = (\pi_\ell,\pi_r) \in S_d \times S_d$ such that the group they generate acts transitively on $\{1,\dots,d\}$. We build two decompositions of each interval $E_i$ as follows. Let $E_{i,\ell} = ]\lambda_{i,\ell},0[$ and $E_{i,r} = ]0,\lambda_{i,r}[$ (which corresponds to the past) and $E_{i,L} = ]\lambda_{i,\ell},\lambda_{i,\ell} + \lambda_{\pi_\ell(i),r}[$, $E_{i,R} = ]\lambda_{i,r}+\lambda_{\pi_r(i),\ell},\lambda_{i,r}[$ (which corresponds to the future). The map $T:E \rightarrow E$ is such that each restriction of $T$ to $E_{i,L}$ (resp. $E_{i,R}$) is a translation onto $E_{\pi_\ell(i),r}$ (resp. $E_{\pi_r(i),\ell}$). We assume implicitly here, that the length-vector $\lambda$ satisfies the \textit{train-track relations}
\begin{equation} \label{eq:train_track_relations}
\lambda_{i,r} - \lambda_{i,\ell} = \lambda_{\pi_\ell(i),r} - \lambda_{\pi_r(i),\ell} \qquad \text{for $i=1,\ldots,d$}.
\end{equation}
We denote by $T_{\pi,\lambda}$ the application constructed above from the data $\lambda$ and $\pi = (\pi_\ell,\pi_r)$ (see the top picture of Figure~\ref{fig:from_train_track_to_L}). We call the map $T: E \rightarrow E$ an \emph{interval exchange transformation}. We warn the reader that our definition of interval exchange transformation does not correspond to the standard one in which one interval is cut in several pieces. In our case many intervals are cut in only two pieces.

Let $T_{\pi,\lambda}:E \rightarrow E$ be an interval exchange transformation. There is a natural way to code the orbits of the dynamical system associated to $T$. To each point $x$ in $E$ we associate a label in $\mathcal{A} = \{1,\ldots,d\} \times \{\ell,r\}$ corresponding to the interval $E_{i,\ell}$ or $E_{i,r}$ it belongs. To an infinite orbit $x$, $T(x)$, $T^2(x)$, \dots we associate an infinite sequence $w \in \mathcal{A}^\NN$ in such way that $w_k$ is the label associated to the point $T^k(x)$. The orbit of $x$ starts with the finite word $w = w_0 w_1 w_2 \ldots w_{N-1} \in \mathcal{A}^N$ if and only if it belongs to the interval
\[
E_w = \bigcap_{k = 0,1,\ldots,N-1} T^{-k} (E_{w_k}).
\]
Similarly to Veech zippered-rectangle construction \cite{Veech82}, we define zippered rectangles for the map $T$.
\begin{definition} \label{def:train_track_suspension}
A \textit{suspension datum} for an interval exchange transformation $(\lambda,\pi)$ is a vector $\zeta = (\zeta_{i,\ell},\zeta_{i,r}) \in (\CC^2)^d$ such that
\begin{itemize}
\item $\Re(\zeta_{i,\ell}) = \lambda_{i,\ell}$ and $\Re(\zeta_{i,r}) = \lambda_{i,r}$,
\item $\Im(\zeta_{i,l}) > 0$ and $\Im(\zeta_{i,r}) > 0$,
\item $\zeta_{i,r} - \zeta_{i,\ell} = \zeta_{\pi_\ell(i),r} - \zeta_{\pi_r(i),\ell}$ (\emph{train-track relations for a suspension}).
\end{itemize}
\end{definition}
To a suspension datum $\zeta$ of $(\pi,\lambda)$ we associate a translation surface in the following way. For $i=1,2,\ldots,d$, let $R_i$ be the quadrilateral with vertices (in trigonometric order) $\zeta_{i,l}$, $0$, $\zeta_{i,r}$, $\zeta_{i,r}+\zeta_{\pi_r(i),\ell} = \zeta_{i,\ell} + \zeta_{\pi_\ell(i),r}$. The suspension $S(\pi,\zeta)$ is the disjoint union of the rectangles $R_i$ for $i=1,\dots,d$ in which we identify the sides which have the same label $(i,\ell)$ or $(i,r)$ (see the second picture in Figure~\ref{fig:from_train_track_to_L}).

\begin{definition}
Let $X$ be a translation surface. A \emph{quadrangulation} of $X$ is a simplicial decomposition of $X$ for which the vertices are the conic singularities of $X$, the edges are geodesics and every face is a quadrilateral which does not contain any singularity. A quadrangulation is \emph{admissible} if every face has exactly two adjacent edges for which the linear flow (in the vertical direction) is incoming.
\end{definition}
In the suspension $S(\pi,\zeta)$ of an interval exchange transformation, the rectangles $R_i$ naturally define an admissible quadrangulation. Reciprocally any admissible quadrangulation gives rise to a train-track: we associate to a quadrangulation the first return map on its sides (see the third picture of Figure~\ref{fig:from_train_track_to_L}).

Let $X$ be a translation surface with an admissible quadrangulation. To an orbit of the linear flow, we associate its natural \emph{cutting sequence} made of the ordered list of edges meet by the orbit. The cutting sequences in a suspension $S(\pi,\zeta)$ of an interval exchange transformation $T_{\pi,\lambda}$ are in bijections with the coding of the orbits in $T_{\pi,\lambda}$.

\begin{figure}[!ht]
\begin{center}\includegraphics{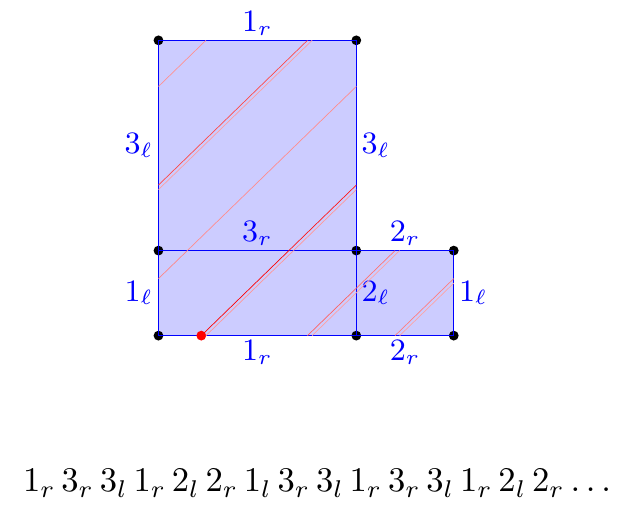}\end{center}
\caption{A geodesic in $\Lsurf(a,b)$ and its cutting sequence.}
\label{fig:Lab_geodesic}
\end{figure}

\subsection{The quadrangulation of $\Lsurf(a,b)$}
We are now interested in the cutting sequences of the linear flow in the surfaces $r_{-\theta} \cdot \Lsurf(a,b)$ defined in Section~\ref{subsection:translation_surfaces} and where $r_{-\theta}$ denotes the rotation by an angle $\theta$. The linear flow in the surface $r_{-\theta} \cdot \Lsurf(a,b)$ is the geodesic flow of $\Lsurf(a,b)$ in the direction $\pi/2 + \theta$.

The next proposition asserts that there is a one to one correspondence between length parameters $((\lambda_{1,\ell},\lambda_{1,r}), (\lambda_{2,\ell},\lambda_{2,r}), (\lambda_{3,\ell},\lambda_{3,r}))$ satisfying the train-track relations for $\pi_\ell = (1,3)$ and $\pi_r=(1,2)$ and parameters $(a,b,\theta)$ of the flat surface $r_{-\theta} \cdot \Lsurf(a,b)$ (up to rescaling).
\begin{proposition} \label{prop:from_L_to_train_track}
Let $\pi_l = (1,2)(3)$, $\pi_r = (1,3)(2)$ and $\lambda$ satisfy the train track relations~(\ref{eq:train_track_relations}). Then there exists a unique suspension datum $\zeta$ of $(\pi,\lambda)$ such that the suspension $S(\pi,\zeta)$ is isomorphic, up to horizontal and vertical rescaling, to a surface of the form $r_\theta \cdot \Lsurf(a,b)$ with $\theta \in ]0,\pi/2[$ endowed with its natural quadrangulation. Moreover, $a$, $b$ and $\theta$ are deduced from the length-vector $\lambda$ and by the following relations
\[
\lambda_{1,\ell} = \lambda_{2,\ell} = -(1-b) \sin \theta  \quad \lambda_{3,\ell} = -b \sin \theta
\qquad
\lambda_{1,r} = \lambda_{3,r} = (1-a) \cos \theta  \quad \lambda_{2,r} = a \cos \theta .
\]
\end{proposition}

\begin{proof}
For $i=1,2,3$, the real parts $\zeta_{i,\ell}$ and $\zeta_{i,r}$ of the suspension datum $\zeta$ are respectively $\lambda_{i,\ell}$ and $\lambda_{i,r}$. We are looking for the imaginary parts such that the suspension $S(\pi,\zeta)$ is isomorphic to a translation surface of the form $r_{\theta} \cdot S(\pi,\zeta)$.

There are two independent equations on the imaginary parts imposed by the train track relations
\begin{equation} \label{eq:tree_relation_in_X}
\zeta_{1,\ell} = \zeta_{2,\ell} \qquad \text{and} \qquad \zeta_{1,r} = \zeta_{3,r}.
\end{equation}
Furthermore, as in any surface $r_\theta \cdot \Lsurf(a,b)$ the sides of the quadrangulation are orthogonal, there are three other independent equations
\begin{equation} \label{eq:perp_relation_in_X}
\zeta_{i,\ell} \perp \zeta_{i,r} \quad \text{for $i=1,2,3$}.
\end{equation}
The equations~(\ref{eq:tree_relation_in_X}) and~(\ref{eq:perp_relation_in_X}) give $5$ independent linear relations for our six parameters $\Im(\zeta_{i,m}), \Im(\zeta_{i,p})$. Hence there is exactly one solution up to rescaling.
\end{proof}

\subsection{$\Lsurf(a,b)$ with extra symmetries: Calta-McMullen L's} \label{subsetion:Veech_surfaces_in_genus2}
Let $X$ be a translation surface with singularities $\Sigma \subset X$. An \emph{affine diffeomorphism} of $X$ is an orientation preserving homeomorphism of $X$ that permutes the singularities of the flat metric and acts affinely on the flat structure of $X$. We denote by $\Aff(X)$ the group of affine diffeomorphisms of $X \backslash \Sigma$. Let denote by $\Gamma(X)$ the image of the derivative map
\[
d : \left\{
\begin{array}{lll}
\Aff(X) & \rightarrow & \SL(2,\RR) \\
f & \mapsto & df
\end{array} \right.
\]
which is called the \emph{Veech group}. For surfaces in $\HHH(2)$, the affine group is isomorphic to the Veech group under the derivative map (see Proposition~4.4 in \cite{HubertLelievre06}). Next, we identify the Veech group and the affine group for surfaces in $\HHH(2)$.

A translation surface $X$ is called a \emph{Veech surface} (or \emph{lattice surface}) if $\Gamma(X)$ is a lattice in $\SL(2,\RR)$. Veech surfaces were introduced in \cite{Veech89} for dynamical purposes. Let $X$ be a Veech surface. For an angle $\theta \in S^1$ if there exists an orbit of the linear flow which joins two singularities (\emph{saddle connection}) then the linear flow $\phi^\theta_t:X \rightarrow X$ in direction $\theta$ is \emph{parabolic}: any geodesic in direction $\theta$ is either a saddle connection or a loop and moreover there exists a non trivial element $\phi \in \Aff(X)$ which stabilizes all geodesics in that direction \cite{Veech89}. The name parabolic comes from the fact that $d \phi \in \SL(2,\RR)$ is a parabolic matrix.

In the surface $\Lsurf(a,b) \in \HHH(2)$ the horizontal and vertical directions are \emph{completely periodic}: all trajectories are either saddle connections or closed loops. A necessary condition for $\Lsurf(a,b)$ to be a Veech surface is that those two directions are parabolic, in other words there exists non identity matrices of the form $\begin{pmatrix}1 & * \\ 0 & 1\end{pmatrix}$ and $\begin{pmatrix}1&0\\ *&1\end{pmatrix}$ in $\Gamma(X)$. It turns out that these conditions are equivalent (which is a miracle of genus $2$ translation surfaces).
\begin{theorem}[\cite{Calta04},\cite{McMullen03}] \label{prop:Veech_L_ab}
The following conditions are equivalent:
\begin{enumerate}
\item the surface $\Lsurf(a,b)$ is a Veech surface,
\item horizontal and vertical directions in $\Lsurf(a,b)$ are parabolic,
\item either $a$ and $b$ are rational or there exists rational numbers $x$ and $y$ and a square-free integer $D > 1$ such that $1/(1-a) = x + y \sqrt{D} $ and $1/(1-b) = (1-x) + y \sqrt{D}$. \label{item:formula}
\end{enumerate}
\end{theorem}
Now, we explicit the form of the parabolic elements in horizontal and vertical directions for parameters that satisfy condition 3 in the above Theorem. Those elements will be used to find paths in the Ferenczi-Zamboni induction. Let $0 < a < 1$ and $0 < b < 1$. Then, the stabilizers in $\SL(2,\RR)$ of, respectively, the bottom and top cylinders in the horizontal decomposition of $\Lsurf(a,b)$ are the two parabolic subgroups of $\SL(2,\RR)$ generated by
\[
P_{bot} = \begin{pmatrix}1 & 1/(1-b)\\0&1\end{pmatrix}
\quad \text{and} \quad
P_{top} = \begin{pmatrix}1&(1-a)/b\\0&1\end{pmatrix}.
\]
Remark that the matrix $P_{bot}$ (resp. $P_{top}$) acts as a Dehn twist around the circumference of the bottom (resp. top) cylinder. The intersection $\langle P_{bot} \rangle \cap \langle P_{top} \rangle$ is non trivial (different from $\{1,-1\}$), if and only if there exist relatively prime positive integers $m_h$ and $n_h$ such that $m_h / (1-b) = n_h (1-a)/b$. The latter equation can be written as $m_h b = n_h (1-a)(1-b)$. By symmetry, the vertical direction is parabolic if and only if there exist relatively prime positive integers $m_v$ and $n_v$ such that $m_v a = n_v (1-a)(1-b)$. We call the 4-tuple $(m_h,n_h,m_v,n_v)$ the \emph{affine multi-twist parameters} of the Veech surface $\Lsurf(a,b)$. It is easy to show that the the existence of $(m_h,n_h,m_v,n_v)$ for parameters $a$ and $b$ is equivalent to the third condition in Proposition~\ref{prop:Veech_L_ab} and more precisely
\begin{proposition}[\cite{Calta04},\cite{McMullen03}]
Let $m_h$, $n_h$, $m_v$ and $n_v$ be positive integers with $m_h$ and $n_h$ (resp. $m_v$ and $n_v$) relatively primes. Then there exist unique real numbers $a$ and $b$ such that $0< a <1$, $0 < b < 1$ and $\Lsurf(a,b)$ is a Veech surface with affine multi-twists parameters $(m_h,m_v,n_h,n_v)$. If we denote $\mu_h = n_h / m_h$ and $\mu_v = n_v / m_v$ then
\begin{align*}
\frac{1}{1-a} = \frac{1 + (\mu_h - \mu_v) + \sqrt{1 + (\mu_h-\mu_v)^2 + 2(\mu_h + \mu_v)}}{2} \\
\frac{1}{1-b} = \frac{1 + (\mu_v - \mu_h) + \sqrt{1 + (\mu_h-\mu_v)^2 + 2(\mu_h + \mu_v)}}{2}
\end{align*}
\end{proposition}
\noindent In particular, $\Lsurf(a,a)$ is a Veech surface if and only if there exists $D = n/m \in \QQ$ with
\[
\frac{1}{1-a} = \frac{1 + \sqrt{1 + 4D}}{2}.
\]
For a parameter $a$ as above, the affine multi-twists parameters $(m_h,m_v,n_h,n_v)$ of $\Lsurf(a,a)$ are $m_h = m_v = m$ and $n_h = n_v = n$.

\subsection{From $\Lsurf(a,b)$ to $T(a,b)$: the wind-tree cocycle} \label{subsection:windtree_cocycle}
Now, we describe a discretization of the billiard flow in $T(a,b)$ as a $\ZZ^2$-cocycle over an interval exchange transformation.

Given a rational billiard, there is a classical procedure to get a translation surface called~\emph{Katok-Zemliakov construction} or \emph{unfolding procedure} (see the original articles~\cite{FoxKershner36} and~\cite{KatokZmeliakov75} or the surveys~\cite{Tabachnikov95} or \cite{MasurTabachnikov02}). The unfolding procedure consists in taking reflected copies of the billiard instead of considering a reflected trajectory. The construction applies to the infinite billiard table $T(a,b)$ and is made of four copies associated to the four directions that a trajectory may take with a given initial angle. We denote by $\Xsurf_\infty(a,b)$ the translation surface obtained by unfolding the billiard table $T(a,b)$.
\begin{proposition}[\cite{DelecroixHubertLelievre11}]
The infinite translation surface $\Xsurf_\infty(a,b)$ is a normal cover of $\Lsurf(a,b)$ and the Deck group $\Deck(\Xsurf_\infty(a,b),\Lsurf(a,b))$ is isomorphic to the semi-direct product $\ZZ^2 \rtimes K$ where $K=\ZZ/2 \times \ZZ/2$ denotes the Klein four group. The intermediate quotient $\Xsurf_\infty(a,b) / \ZZ^2$ is a four fold cover of $\Lsurf(a,b)$ which corresponds to the unfolding of the billiard in a fundamental domain of $T(a,b)$ under the $\ZZ^2$-action.
\end{proposition}
The proof of the proposition is elementary and we refer to~\cite{DelecroixHubertLelievre11}. Now, we consider how a geodesic in $\Xsurf_\infty(a,b)$ can be built from the ones in $\Lsurf(a,b)$. The Klein four group $K = \ZZ/2 \times \ZZ/2$ in the proposition naturally identifies with the group generated by the vertical and horizontal reflexions denoted respectively $\tau_v$ and $\tau_h$.

\begin{figure}[!ht]
\begin{center}\includegraphics{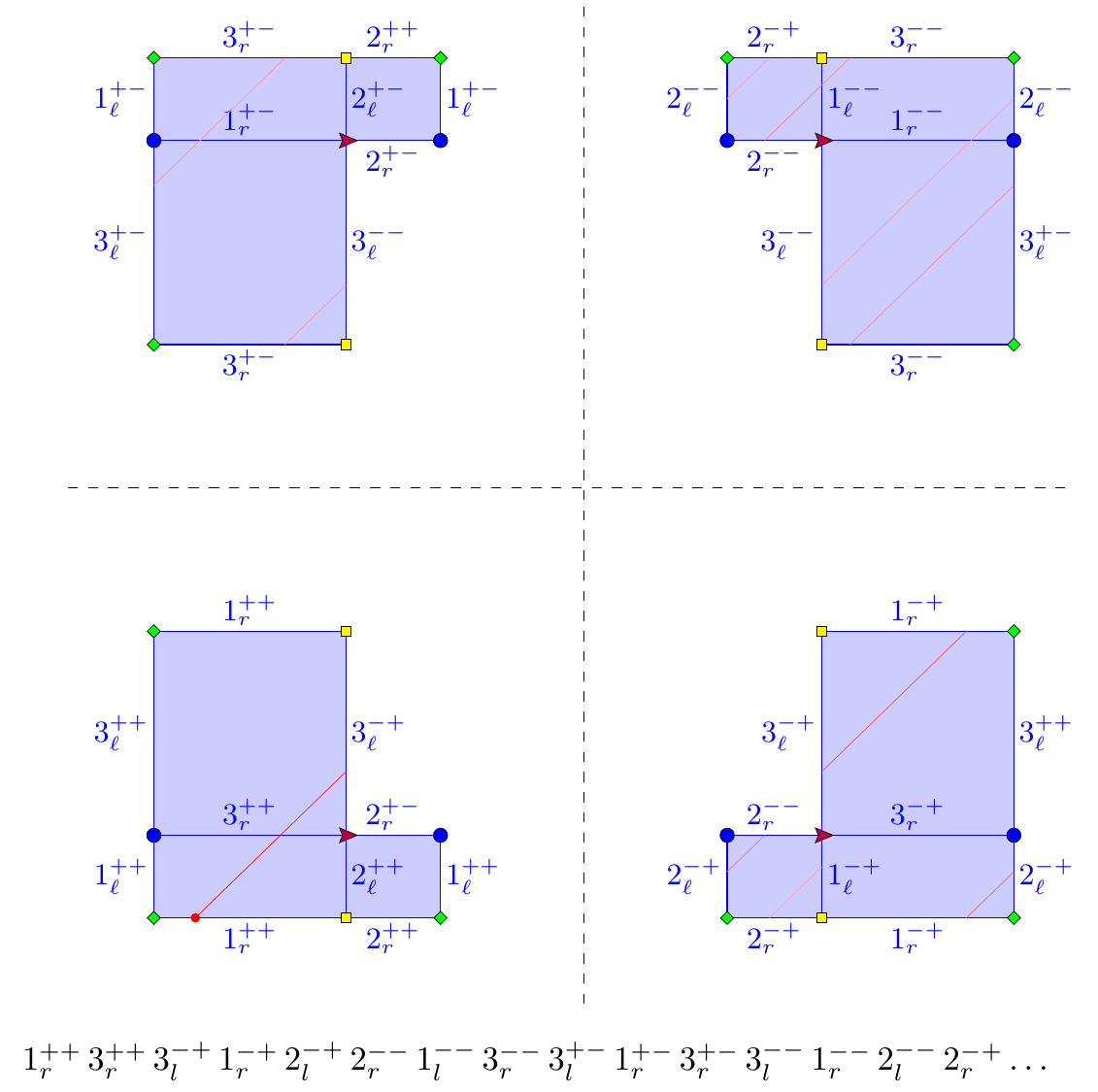}\end{center}
\caption{The lift in $\Xsurf(a,b)$ of the geodesic in $\Lsurf(a,b)$ from Figure~\ref{fig:Lab_geodesic}.}
\label{fig:Xab_geodesic}
\end{figure}

In Section~\ref{section:train_tracks_and_quadrangulations}, we defined a symbolic coding of geodesics in $\Lsurf(a,b)$ on the alphabet $\mathcal{A} = \{1,2,3\} \times \{\ell,r\}$. The preimage of the quadrangulation of $\Lsurf(a,b)$ in $\Xsurf(a,b)$ gives a symbolic coding on the alphabet $\mathcal{A} \times K$. Given a geodesic $\gamma$ (finite or infinite) in $\Lsurf(a,b)$ and its coding $w = (w_i)_i$ (in $\mathcal{A}^*$ or $\mathcal{A}^\NN$) it has four lifts in $\Xsurf(a,b)$ and hence four possible codings. The \emph{canonical} one that we denote $\tilde{w} = ((w_i,\kappa_i))_i$ is the one which starts in the copy labelled $id \in K$. Let $g: \mathcal{A} \rightarrow K$ be defined by
\[
g(1_\ell) = g(1_r) = g(2_\ell) = g(3_r) = id \qquad g(2_r) = \tau_h \qquad \text{and} \qquad g(3_\ell) = \tau_v.
\]
Then the canonical lift $\tilde{w} = ((w_i,\kappa_i))$ of $w$ can be defined by
\[
\kappa_0 = 1 \quad \text{and for $i \geq 0$, } \quad \kappa_{i+1} = \kappa_i\, g(w_i).
\]
The Klein four group $K$ acts transitively on the four lifts of $w$. The three other cutting sequences for the lifts are $((a_i, \tau_v \kappa_i))_i$, $((a_i, \tau_h \kappa_i))_i$ and $((a_i, \tau_v\,\tau_h \kappa_i))_i$. In Figure~\ref{fig:Xab_geodesic} we give an example of the lift of a geodesic in $\Lsurf(a,b)$. To simplify notations we use $++$ (resp. $-+$, $+-$ and $--$) for the element $id \in K$ (resp. $\tau_v$, $\tau_h$ and $\tau_v \tau_h$).

The surface $\Xsurf_\infty(a,b)$ is a $\ZZ^2$-cover of $\Xsurf(a,b)$. Hence, the preimage of the quadrangulation of $\Xsurf(a,b)$ determines a quadrangulation in $\Xsurf_\infty(a,b)$. We fix an origin in $\Xsurf_\infty(a,b)$ and consider a bijection between the faces of the quadrangulation and $\{1,2,3\} \times K \times \ZZ^2$. To lift the cutting sequence of a geodesic in $\Xsurf(a,b)$ to $\Xsurf_\infty(a,b)$ there is another cocycle which is defined on the copies $\{1,2,3\} \times \{id\}$ by
\[
\begin{array}{c}
f((2_\ell,id)) = f((2_r,id)) = f((3_\ell,id)) = f((3_r,id)) = (0,0)\\
f((1_\ell,id)) = (1,0) \qquad \text{and} \qquad f((1_r,id)) = (0,1)
\end{array}
\]
and on the three other copies by the symmetry rule
\[
\forall a \in \{1,2,3\}, \forall \kappa \in K,\ f((a,\kappa)) = \kappa \cdot f((a,id))
\]
where $\tau_h$ and $\tau_v$ acts on $\ZZ^2$ by reflexion
\[
\tau_h \cdot (x,y) = (x,-y) \qquad \tau_v \cdot (x,y) = (-x,y).
\]

\begin{proposition}[\cite{DelecroixHubertLelievre11}]
Let $\gamma$ be a geodesic in $\Xsurf_\infty(a,b)$, $w$ its cutting sequence on $\{1,2,3\} \times K \times \ZZ^2$ and $\overline{\gamma}$ its image in $T(a,b)$. Then for all $t > 0$ we have
\[
\left| d(x,\phi_t^\theta(x)) - \|f^{(n)}(w)\|_2  \right| \leq \sqrt{2}
\]
where $f^{(n)}(w) = f(w_0) f(w_1) \ldots f(w_{n-1})$ and $n$ is such that the geodesic from $x$ has cut $n$ sides of the quadrangulation before time $t$.
\end{proposition}

As the cover $\Xsurf_\infty(a,b) \rightarrow \Lsurf(a,b)$ is normal, we can build a non-commutative cocycle to lift the cutting sequences of geodesics in $\Xsurf(a,b)$ to $\Xsurf_\infty(a,b)$.

To simplify notations, we use a direct description from $\Lsurf(a,b)$ to $\Xsurf_\infty(a,b)$. Let $D_\infty = \ZZ \rtimes \ZZ/2$ be the infinite dihedral group (where $\ZZ/2$ acts by multiplication by $-1$ on $\ZZ$). We use the following notation for $G = D_\infty^2$. The generators of $K=\ZZ/2 \times \ZZ/2$ are denoted by $\tau_v$ and $\tau_h$. We use multiplicative notations and write the product rule in $D_\infty^2$ as follows. For $(x_1,y_1),(x_2,y_2),(x,y) \in \ZZ^2$ we have
\[
(x_1,y_1)(x_2,y_2) = (x_1+x_2,y_1+y_2) \qquad \tau_v (x,y) = (-x,y) \tau_v \qquad \tau_h (x,y) = (x,-y) \tau_h.
\]
The cocycle which describe a cutting sequence in $\Xsurf_\infty(a,b)$ from one in $\Lsurf(a,b)$ is the map $f:\mathcal{A} \rightarrow G$ defined by
\[
\begin{array}{lll}
f(1_\ell) = (1,0) & f(2_\ell) = \tau_v & f(3_\ell) = (0,0) \\
f(1_r) = (0,1) & f(2_r) = (0,0) & f(3_r) = \tau_h.
\end{array}
\]
The cocycle $f$ can be viewed as a function $\phi$ on the domain $E = E_1 \sqcup E_2 \sqcup E_3$ of an interval exchange transformation $T = T_{(\lambda,\pi)}$ with $\pi_l=(1,3)$ and $\pi_r=(1,2)$ which is constant on each interval $E_{i,\ell}$ (resp. $E_{i,r}$). Its value on $E_{i,\ell}$ is given by $f(i_\ell)$ (resp. by $f(i_r)$).

\section{Divergent wind-tree cocycles} \label{section:divergent_cocycle}
We recall the renormalization procedure of S.~Ferenczi and L.~Zamboni~\cite{FerencziZamboni10_eig,FerencziZamboni10_struct} for linear flow on translation surfaces. The induction was used to obtain fine properties of interval exchanges in the hyperelliptic classes and many examples of exotic ergodic behaviors of interval exchange transformations. We use their induction to control the wind-tree cocycle. An other induction procedure for interval exchange transformations is the one of G.~Rauzy~\cite{Rauzy79} which seems less adapted to our situation.

\subsection{Ferenczi-Zamboni induction in hyperelliptic strata} \label{subsection:FZ_induction}
We now recall the induction procedure introduced in~\cite{FerencziZamboni10_struct}. In next sections, we restrict our study to the case of the stratum $\HHH(2)$ which is the subject of~\cite{FerencziZamboni10_eig} and corresponds to our surface $\Lsurf(a,b)$.

Let $X$ be a translation surface with an admissible quadrangulation. The general principle of the Ferenczi-Zamboni induction consists in looking at a sequence of admissible quadrangulations of the surface such that the quadrilaterals become more and more flat in the direction of the linear flow. In the case of hyperelliptic strata, we consider only quadrangulations which are stable under the hyperelliptic involution. This restriction guarantees the existence of an induction procedure. The main point is that the hyperelliptic involution simplifies the train-track relations~(\ref{eq:train_track_relations}) (see Definition~2.4 and the discussion which follows in \cite{FerencziZamboni10_struct}).

\begin{figure}[!ht]
\begin{center}
\picinput{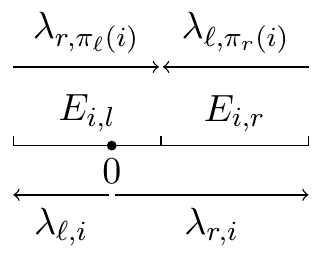} \hspace{2cm} \picinput{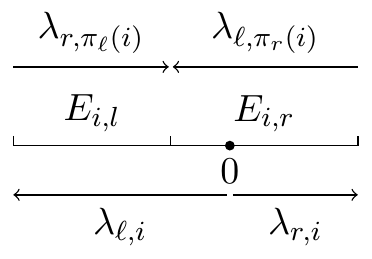}
\end{center}
\caption{Left or right state for the interval $E_i$ of an interval exchange transformation.}
\label{fig:left_or_right_induction}
\end{figure}

Now, we describe the induction. Let $T=T_{\pi,\lambda}:E \rightarrow E$ be an interval exchange transformation on $d$ intervals. We assume that $E = E_1 \sqcup \ldots \sqcup E_d$ is stable under the hyperelliptic involution. We want to define a new interval exchange transformation which corresponds to a first return map on a union of $d$ subintervals $E' = E'_1 \sqcup \ldots \sqcup E'_d \subset E$ where for each $i=1,\ldots,d$, $E'_i \subset E_i$ and $E'_i$ contains the origin of $E_i$. For each $i$, we define its \emph{state} (see Figure~\ref{fig:left_or_right_induction}):
\begin{itemize}
\item $i$ is in \emph{left state} if $\lambda_{\pi_\ell(i),r} > \lambda_{i,\ell}$ (or equivalently $0 \in E_{i,l}$),
\item $i$ is in \emph{right state} if $\lambda_{\pi_r(i),\ell} > \lambda_{i,r}$ (or equivalently $0 \in E_{i,r}$).
\end{itemize}
\noindent Knowing the state of each level we want to define $T'$ by choosing among the following choices
\begin{itemize}
\item if $i$ is in a left state either we choose $E'_i = E_i$ or $E'_i = E_{i,\ell}$,
\item if $i$ is in a right state either we choose $E'_i = E_i$ or $E'_i = E_{i,r}$.
\end{itemize}

Now, let $T=T_{\pi,\lambda}$ be an interval exchange transformation in a hyperelliptic strata such that it is stable under the hyperelliptic involution. As shown in~\cite{FerencziZamboni10_struct} there are canonical choices which ensure that first of all the induction is always possible. Secondly, the induced interval exchange transformation is also stable under the hyperelliptic involution. The choice is done as follows. We consider the cycles in the disjoint cycle decomposition of the permutations $\pi_\ell$ and $\pi_r$. If there exists a cycle $c$ of $\pi_\ell$ for which each element of $c$ are in left state then we allow to perform a left induction for all of them. Such a cycle is called a \emph{left branch of induction}. Formally a left induction step for $c$ on $\pi = (\pi_\ell,\pi_r)$ and $\lambda = (\ell,r)$ gives the combinatorial data $\pi' = (\pi_\ell,\pi'_r)$ and $\lambda' = (\lambda_\ell,\lambda'_r)$ where
\begin{itemize}
\item for all $i$ in $c$, $\lambda'_{i,r} = \lambda_{i,r} + \lambda_{\pi_r(i),\ell}$ and $\pi'_r(i) = \pi_r \circ \pi_\ell(i)$,
\item for all $i$ not in $c$, $\lambda'_{i,r} = \lambda_{i,r}$ and $\pi'_r(i) = \pi_r(i)$.
\end{itemize}
The definition of induction can easily be extended to suspensions by replacing $\lambda_{i,\ell}$ and $\lambda_{i,r}$ in the formulas above by respectively $\zeta_{i,\ell}$ and $\zeta_{i,r}$.

The following general theorem can be checked by hand for $\HHH(2)$ by building the so called ``graph of graphs''(see Figure~\ref{fig:graph_of_graphs}).
\begin{theorem}[\cite{FerencziZamboni10_struct} Lemma~2.5 and Proposition~2.6] \label{prop:paths_in_graph_of_graphs}
Let $T$ be a hyperelliptic interval exchange in $\HHH^{hyp}(2g-2)$ or $\HHH^{hyp}(g-1,g-1)$ without saddle connection. Then, $T$ admits at least one induction branch and any map induced from $T$ by cutting all intervals in an induction branch is an hyperelliptic interval exchange transformation. Moreover, for any choice of infinite sequence of inductions for $T$ such that
\begin{itemize}
\item if $i$ is not in the induction branch at stage $n$ then the state of $i$ at step $n$ and $n+1$ is the same
\item each interval $E_i$ for $i=1,2,\ldots,d$ is cut infinitely many times on its left and on its right.
\end{itemize}

Conversely, given an infinite path of inductions in the graph of graphs starting from $\pi$ that satisfies the two conditions above, there exists at least one parameter $\lambda$ for which the interval exchange $T_{\pi,\lambda}$ has no saddle connection and from which we can perform these steps of induction.
\end{theorem}
We use the following multiplicative algorithm similar to the Gauss map for coding geodesics in the torus.
\begin{definition} \label{def:multiplicative_algorithm}
The \emph{multiplicative Ferenczi-Zamboni} algorithm on a symmetric interval exchange transformations is the algorithm which at odd steps performs all possible right inductions and even steps all possible left inductions.
\end{definition}

\subsection{Description of the language in terms of induction} \label{subsection:language}
We now follow~\cite{FerencziZamboni10_eig} to describe the language of an interval exchange transformation in terms of one of its induction. Let $T = T_{\pi,\lambda}$ be a symmetric interval exchange transformation in $\HHH(2)$. For $i=1,2,3$, we note $L_i = (\pi_l(i),r)$ and $R_i = (\pi_r(i),\ell)$ and $L=(L_1,L_2,L_3)$ and $R=(R_1,R_2,R_3)$. The words $L_i$ and $R_i$ are the two possible continuations of a letter of the form $(i,*)$.

Let $c$ be a union of left (or a union of right) admissible branch of induction for $T$ and let $T'$ be the interval exchange transformation obtained from $T$ by inducing with respect to $c$. Then the words $L'_i$ and $R'_i$ on this new interval exchange transformation seen from $T$ can be described by the following rules
\begin{center}\begin{tabular}{l|l}
left induction on $c$ & right induction on $c$ \\
\hline
$L'=L$  & $R'=R$\\
for $i$ in $c$, $R'_i = L_i R_{\pi_\ell(i)}$ & for $i$ in $c$, $L'_i = R_i L_{\pi_r(i)}$ \\
for $i$ not in $c$, $R'_i = R_i$ & for $i$ not in $c$, $L'_i = L_i$\\
\end{tabular}\end{center}

\begin{figure}[!ht]
\picinput{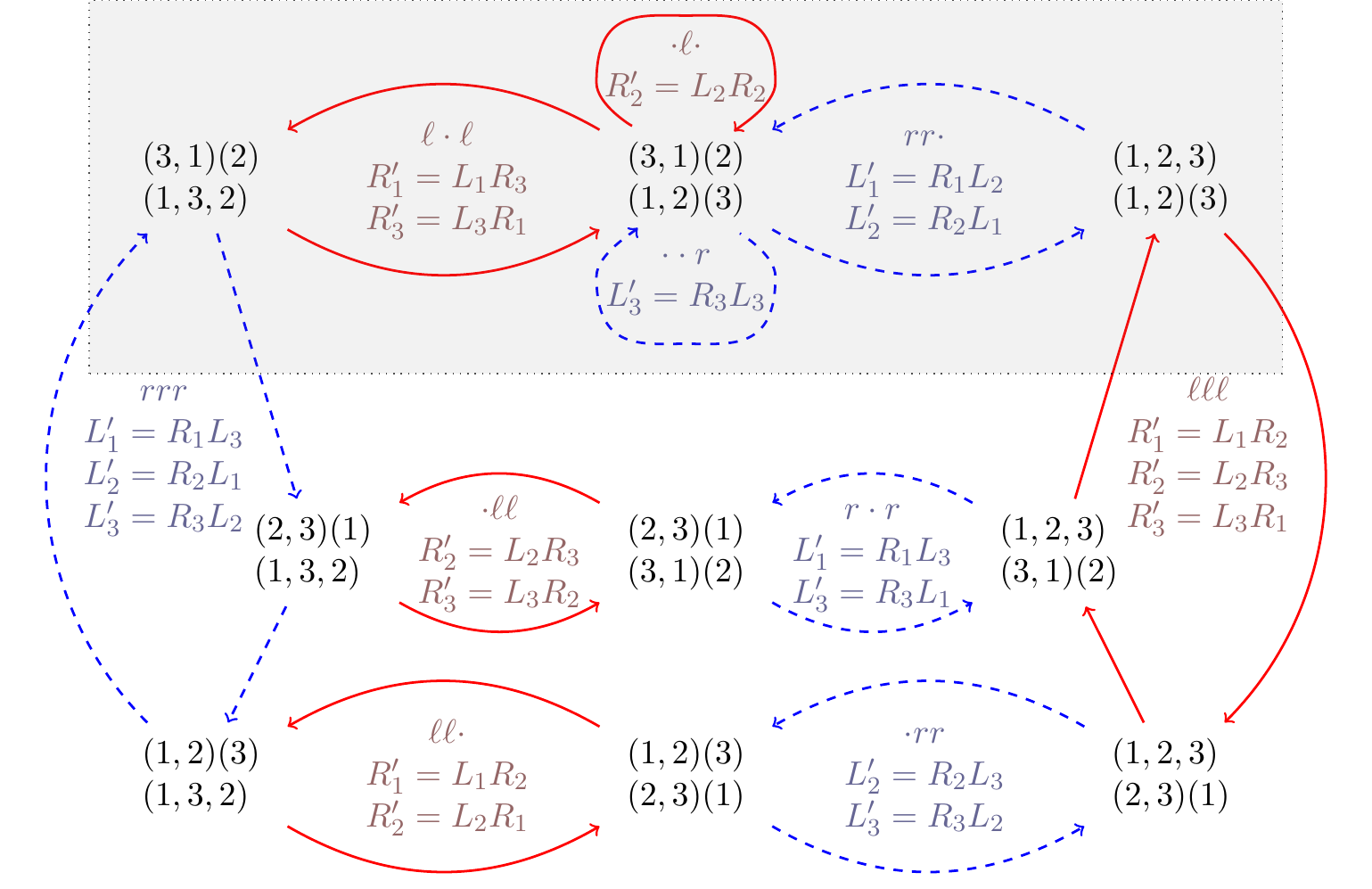}
\caption{The states of the induction procedure in $\HHH(2)$ (left inductions are dashed). Some inductions which correspond to loops are missing. The light colored rectangle corresponds to the subset of inductions considered in Section~\ref{section:subset_of_inductions}.}
\label{fig:graph_of_graphs}
\end{figure}

Starting from a symmetric interval exchange transformation and performing successively left and right inductions, we build a sequence of words $L_i^{(k)}$ and $R_i^{(k)}$. The possible inductions at each step are shown in Figure~\ref{fig:graph_of_graphs}. From the definition of our multiplicative algorithm, at odd step $2k+1$ for $k \geq 0$, we perform a right induction and we have $R_i^{(2k+1)} = R_i^{(2k)}$. Similarly at even step $2k+k$ for $k \geq 0$, we perform left inductions and we have $L_i^{(2k+2)} = L_i^{(2k+1)}$.

\subsection{A subset of parameters $(a,b,\theta)$ defined from the induction} \label{section:subset_of_inductions}
We restrict our attention to a subset of inductions which simplifies considerably the form of the infinite words we obtain. This subset of possible inductions is similar to the one used in~\cite{FerencziZamboni10_eig} Section~5. For parameters $(a,b,\theta)$ associated to these inductions, we will be able to control the billiard orbits of $T(a,b)$ in direction $\theta$.

Our graph of induction consists of the unique state $\pi = (\pi_\ell,\pi_r)$ with $\pi_\ell=(1,3)$ and $\pi_r=(1,2)$ which corresponds to the quadrangulation of $\Lsurf(a,b)$. We consider as induction steps
\begin{itemize}
\item the induction which are a succession of left inductions only that go from $\pi$ to $\pi$,
\item the induction which are a succession of right inductions only that go from $\pi$ to $\pi$.
\end{itemize}
The subgraph of inductions is the part of the graph of graphs in the rectangle in Figure~\ref{fig:graph_of_graphs}. There are two possible left induction from $\pi$ associated respectively to the states $\ell \cdot \ell$ and $\cdot \ell \cdot$. As two left inductions commute, each step of the multiplicative algorithm (Definition~\ref{def:multiplicative_algorithm}) corresponds to a $2$-tuple of integers $(m,n)$ where $m$ corresponds to the multiplicity of the loop of length two associated to state $\ell \cdot \ell$ and $n$ is the multiplicity of the loop associated to the state $\cdot \ell \cdot$. A $2$-tuple also encodes the left inductions. The induction algorithm is then a shift on sequences of two-tuples of integers $((m_k,n_k))_{k \geq 0}$.

The language of an interval exchange transformation that admits an induction which belongs to the subgraph is defined by two families of substitutions indexed by a $2$-tuple $(m,n)$ of integers
\begin{equation} \label{eq:sigma_l_and_sigma_r}
\sigma_\ell(m,n) = 
\left\{ \begin{array}{l@{\ \mapsto \ }l}
L_1 & (R_1 R_2)^m L_1 \\
L_2 & (R_2 R_1)^m L_2 \\
L_3 & (R_3)^n L_3
\end{array} \right.
\qquad \text{and} \qquad
\sigma_r(m,n) =
\left\{ \begin{array}{l@{\ \mapsto \ }l}
R_1 & (L_1 L_3)^m R_1 \\
R_2 & (L_2)^n R_2 \\
R_3 & (L_3 L_1)^m R_3
\end{array} \right.
\end{equation}

Because of the train track relations for $\pi$, namely $\lambda_{1,\ell} = \lambda_{2,\ell}$ and $\lambda_{1,r} = \lambda_{3,r}$, the possible vector-lengths are described by the $3$-tuple $Z = (x_1,x_2,y_1,y_2) = (|\lambda_{2,\ell}|,|\lambda_{3,\ell}|,|\lambda_{3,r}|,|\lambda_{2,r}|) \in \PP^3(\RR)$. The application of one step of the above algorithm corresponds, at the level of vector-lengths, to the application of one of the two projective maps below
\[
F_\ell(m,n): \begin{pmatrix}x_1\\x_2\\y_1\\y_2\end{pmatrix} \mapsto
\begin{pmatrix}
x_1 - m (y_1+y_2) \\
x_2 - n y_1 \\
y_1 \\
y_2
\end{pmatrix}
\qquad
F_r(m,n): \begin{pmatrix}x_1\\x_2\\y_1\\y_2\end{pmatrix} \mapsto
\begin{pmatrix}
x_1 \\
x_2 \\
y_1 - m (x_1+x_2) \\
y_2 - n x_1
\end{pmatrix}.
\]
The multiplicative induction algorithm on the subgraph corresponds exactly to the map $F$ defined in the introduction. We emphasize that not all vector-lengths parameters $Z=(\lambda_{2,\ell},\lambda_{3,\ell},\lambda_{3,r},\lambda_{2,r})$ admit a continued fraction expansion with respect to this algorithm (the domain is a Cantor set). We recall that we name $F$-\textit{renormalizable} a quadruple of length parameters $Z$ for which the induction is exactly prescribed by our subgraph with one vertex at $\pi = ((1,3),(1,2))$. From Theorem~\ref{prop:paths_in_graph_of_graphs}, an $F$-renormalizable $3$-tuple determines a unique sequence of $2$-tuples $((m_k,n_k))_{k \geq 0}$ such that
\begin{itemize}
\item for each $k \geq 1$ either $m_k \not= 0$ or $n_k \not= 0$,
\item if $m_k=0$, then $m_{k+1} \not= 0$ and $n_{k+1} = 0$
\item for infinitely many $i$, $m_{2k} \not = 0$ (resp. $m_{2k+1} \not= 0$),
\item for infinitely many $i$, $n_{2k} \not = 0$ (resp. $n_{2k+1} \not= 0$).
\end{itemize}
Reciprocally, from Proposition~\ref{prop:paths_in_graph_of_graphs}, we know that every sequence of $2$-tuples of integers that satisfy the above conditions gives a $F$-renormalizable vector-lengths.

\subsection{Renormalizable slopes in Veech $\Lsurf(a,b)$} \label{subsection:veech_parameters}
In this section, to a Veech surface of the form $\Lsurf(a,b)$ (see Proposition~\ref{prop:Veech_L_ab} for a characterization) we build a set of slopes $\Lambda \subset S^1$ for which the corresponding vector lengths of the interval exchange transformation $F$-renormalizable.

\begin{proposition} \label{prop:infinitely_renormalizable_veech}
Let $\Lsurf(a,b)$ be a Veech surface and $(m_h,n_h,m_v,n_v)$ its Dehn multi-twist parameters (see Section~\ref{subsetion:Veech_surfaces_in_genus2}). Let $(s_h,s_v)$ be the widths of the associated parabolic matrices in the Veech group. Then for $\theta \in (0,\pi/2)$ of the form
\[
\tan(\theta)^{-1} = a_0\,s_h + \cfrac{1}{a_1\,s_v + \cfrac{1}{a_2\,s_h + \cfrac{1}{a_3\,s_v + \cfrac{1}{\ldots}}}}
\]
the interval exchange transformation $T$ associated to $(a,b,\theta)$ by Proposition~\ref{prop:from_L_to_train_track} is $F$-renormalizable. The convergents associated to the restricted multiplicative Ferenczi-Zamboni induction of $T$ are $(a_{2k} m_h, a_{2k} n_h)$ for even $k$ and $(a_{2k+1} m_v, a_{2k+1} n_v)$ for odd $k$.
\end{proposition}

\begin{proof}
It is more convenient to consider coordinates in $\PP^1(\RR)$ instead of $\theta \in (0,\pi/2)$. To the angle $\theta$ we associate the (oriented) slope $x = \tan(\theta)^{-1} = (\cos(\theta):\sin(\theta)) \in \PP^1(\RR)$.

Let $\rho_\ell$ (resp. $\rho_r$) be the vertical (resp. the horizontal) parabolic element which stabilizes the Veech surface $\Lsurf(a,b)$. We note
\[
\rho_\ell = \begin{pmatrix}1&0\\s_v&1\end{pmatrix}
\qquad \text{and} \qquad
\rho_r = \begin{pmatrix}1&s_h\\0&1\end{pmatrix}.
\]
As the Veech group of $\Lsurf(a,b)$ is a lattice $s_h s_v \geq 1$, otherwise the group generated by $\rho_\ell$ and $\rho_r$ won't be discrete. In particular, if $x=x(\theta)$ has the form given in the statement, then the sequence $(a_k)_k$ associated to $x \in \PP^1(\RR)$ is unique. More precisely, the expansion of $x$ is defined from a modified continued fraction algorithm. Let $\psi_\ell:\left]0,1/s_v\right[ \rightarrow ]1/s_v,\infty[$ and $\psi_r: \left]s_h,\infty \right[ \rightarrow \left]0,s_h \right[$ be the two maps
\[
\psi_\ell(x) = \frac{1}{1/x - \left\lfloor \frac{1/x}{s_v} \right\rfloor s_v}
\quad \text{and} \quad
\psi_r(x) = x - \left\lfloor \frac{x}{s_h} \right\rfloor s_h
\]
As $s_h s_v \geq 1$ the domain of $\psi_\ell$ and $\psi_r$ are disjoint and we define $\psi$ to be the map that equals $\psi_\ell$ on $(0,1/s_v)$ and $\psi_r$ on $(s_h,\infty)$. The map $\psi$ is associated to the shift on the sequence $(a_k)_k$ that defines $x=\tan(\theta)$: if $x$ is defined by the sequence $(a_0,a_1,a_2,\ldots)$ then $\psi(x)$ is defined either by the sequence $(0,a_1,a_2,\ldots)$ if $a_0\not=0$ or $(a_2,a_3,\ldots)$ if $a_0 = 0$.

The maps $\psi_\ell$ and $\psi_r$ correspond to the standard projective action of powers of the inverses of two matrices $\rho_\ell$ and $\rho_r$ that corresponds to the horizontal and vertical multi-twist in $\Lsurf(a,b)$:

Let $\theta \in S^1$ be such that $x=\tan(\theta)^{-1}$ admits an infinite expansion with respect to $\psi$ and $S(\pi,\zeta)$ be the suspension associated to $r_{\pi/2-\theta} \cdot \Lsurf(a,b)$ as in Proposition~\ref{prop:from_L_to_train_track}. Assume that $a_0 \not= 0$, then we can perform a left induction with parameters $(a_0 m_h, a_0 n_h)$. Let $S(\pi,\zeta')$ be the surface obtained after this step of left induction. Then the vertical direction in $S(\pi,\zeta')$ corresponds to the direction $x$ in $\rho_h^{a_0} \cdot \Lsurf(a,b)$ or equivalently, to the direction $\psi(x)$ in $\Lsurf(a,b)$.
\end{proof}

We now estimate the Hausdorff dimension of the set of renormalizable slope in a Veech surface of the form $\Lsurf(a,b)$.
\begin{proposition}
Let $s$ and $t$ be positive real numbers with $s t \geq 1$ and consider $R(s,t) \subset \RR$ the set of real numbers $x$ of the form
\[
x = a_0\,s + \cfrac{1}{a_1\,t + \cfrac{1}{a_2\,s + \cfrac{1}{a_3\,t + \cfrac{1}{\ldots}}}}
\]
where $a_0 \geq 0$ and $a_k \geq 1$ for $k \geq 1$ are integers. Then
\begin{itemize}
\item for any $\lambda > 0$, $\Hdim(R(s_h,s_v)) = \Hdim(R(\lambda s, \lambda^{-1} t))$,
\item the map $s \mapsto \Hdim(R(s,s))$ is decreasing and not smaller than $1/2$.
\end{itemize}
\end{proposition}

\begin{proof}
We use the following notation
\[
[a_0,a_1,\ldots]_{s,t} = a_0\,s + \cfrac{1}{a_1\,t + \cfrac{1}{a_2\,s + \cfrac{1}{a_3\,t + \cfrac{1}{\ldots}}}}.
\]
The map
\[
\begin{array}{ccc}
R(s,t) & \rightarrow & R(\lambda s, \lambda^{-1} t) \\
x=[a_0,a_1,\ldots]_{s,t} & \mapsto & [a_0,a_1,\ldots]_{\lambda s, \lambda^{-1} t}
\end{array}
\]
is just a multiplication by $\lambda$. As bi-Lipschitz map preserves Hausdorff dimension, $\Hdim R(s,t) = \Hdim R(\lambda s, \lambda^{-1} t)$.

The fact that $s \mapsto \Hdim R(s,s)$ is decreasing is immediate from the construction. In the following, we fix $s > 1$ and show that the Hausdorff dimension of $R(s,s)$ is not smaller than $1/2$. It follows from \cite{Hensley} Chapter~9 that the Hausdorff dimension can be computed with the canonical covers. More precisely, for a tuple $v=(a_1,\ldots,a_k)$ we define the $s$-convergents as follows
\[
\begin{array}{lll}
p_0 = 1  & p_1 = s\,a_2 & p_{k+1} = s\, a_k p_k + p_{k-1} \\
q_0 = s\,a_1 & q_1 = s^2\,a_2 a_1 + 1 & q_{k+1} = s\, a_k q_k + q_{k-1}.
\end{array}
\]
Then the quantity $[a_1,a_2\ldots,a_k]_{s,s} = p_k/q_k$ converges to the real number
\[
x = \cfrac{1}{a_1\,s + \cfrac{1}{a_2\,s + \cfrac{1}{a_3\,s + \cfrac{1}{\ldots}}}}.
\]
Let $\mathbf{v}=(a_1,\ldots,a_k)$ be a $k$-tuple of integers. We denote by $|\mathbf{v}|_s$ the denominator $q_k$ of the continued fraction $[s\,a_1,\ldots,s\,a_k]$ constructed above. To define the Hausdorff dimension we first define the function $\lambda$ as follows
\[
\lambda_s(\sigma) = \sup \left\{\lambda > 0; \limsup_{r \to \infty} \lambda^{-r} \sum_{\mathbf{v} \in \NN^r} |\mathbf{v}|_s^{-\sigma} \right\}
\]
where $\NN$ is the set of positive integers. The length of the interval $\{[a_1,\ldots,a_k+t]_s; t \in [0,1]\}$ is $|(a_1,\ldots,a_k)|_s^{-1} |(a_1,\ldots,a_k+1)|_s^{-1}$. Hence, the quantity $|\mathbf{v}|_s^{-\sigma}$ in the definition of $\lambda$ is up to a factor $2$, the length of an interval of the canonical cover at step $r$ to the power $\sigma/2$. The Hausdorff dimension of $R(s,s)$ is
\[
\Hdim R(s,s) = \frac{1}{2} \inf \{\sigma; \lambda_s(\sigma) < 1\}.
\]
But $\lambda_s(1) = 0$ as the serie $\displaystyle \sum_{\mathbf{v} \in \NN^r} |\mathbf{v}|_s^{-1}$ diverges for any $s$. The Hausdorff dimension of $R(s,s)$ is then not smaller than $1/2$.
\end{proof}

\begin{corollary}
For any Veech $\Lsurf(a,b)$ with parabolics $\rho_\ell$ and $\rho_r$ in respectively vertical and horizontal direction, the set of slopes renormalizable by $\rho_\ell^m$ and $\rho_r^n$ has positive Hausdorff measure bounded below by $1/2$ for any positive integer $m$ and $n$.
\end{corollary}

\subsection{Divergent cocycles: proof of the main theorem} \label{subsection:main_proof}
This section is devoted to the proof of Theorem~\ref{thm:even_implies_div}. In next section, we illustrate all computations with the simple example of $L(1/2,1/2)$ with slope $\theta=\arctan(\sqrt{2}-1)$ which corresponds to the periodic expansion $((1,2),(1,2),\ldots)$.

We first describe the strategy of the proof. Let $T$ be an interval exchange transformation without saddle connection that is $F$-renormalizable. Then, for each step $k \in \NN$, the $k$-th step of the Ferenczi-Zamboni induction can be used to decompose the coding of a trajectory with the six pieces $L_i^{(k)}$ and $R_i^{(k)}$ for $i=1,2,3$. The size of the pieces grows with $k$ and more precisely, the pieces at step $k+1$ are concatenations of the pieces at step $k$. The rule to glue the pieces is given by the substitutions defined in Section~\ref{subsection:language} and depends on the convergents $((m_k,n_k))_{k \geq 0}$ of the restricted Ferenczi-Zamboni induction of $T$. Now assume that $T$ satisfies the statement of Theorem~\ref{thm:even_implies_div}. Then we prove that for each $k \geq 0$ the wind-tree cocycle over $T$ has no ``local self-intersection''. More precisely, for $i=1,2,3$ and $k \geq 0$, let $\LLL_i^{(k)}$ and $\RRR_i^{(k)}$ be the subsets of $G = D_\infty \times D_\infty$ made of the values taken by the wind-tree cocycle on the finite pieces $L_i^{(k)}$ and $R_i^{(k)}$. In the cutting sequence considered as concatenation of pieces $L_i^{(k)}$ and $R_i^{(k)}$, the values taken by the cocycles on all pieces are translates of $\LLL_i^{(k)}$ and $\RRR_i^{(k)}$ by an element of $G = (\ZZ \rtimes \ZZ/2)^2$. We prove that for $k > 2$, the values of level $k+1$ are built in such way that each part from level $k$ do not intersect each other. The reason why we need $k > 2$ is due to the fact that for step $1$ (resp. step $2$) the trajectory can rebound between two vertical scatterers (resp. horizontal scatterers) which implies that the values of the cocycle during this period take only two values $(x,y)$ and $(x,y)\tau_v$ (resp. $(x,y)$ and $(x,y)\tau_h$). In particular we prove the stronger statement that the trajectory in the wind-tree model are ``self-avoiding''.

We fix for the remaining of the section a triple $(a,b,\theta)$ that fulfill the hypothesis of Theorem~\ref{thm:even_implies_div} and consider the associated interval exchange transformation $T$. We denote by $((m_k,n_k))_{k \geq 0}$ the convergents of the restricted Ferenczi-Zamboni induction of $T$ and $(L^{(k)},R^{(k)})$ the $6$-tuple of words that describe the coding of the orbits in $T$. For $i=1,2,3$, let $\LLL_i^{(k)}$ and $\RRR_i^{(k)}$ the subsets of $G = D_\infty \times D_\infty$ made of values taken by the wind-tree cocycle on respectively $L_i^{(k)}$ and $R_i^{(k)}$. At step $k=0$
\begin{equation} \label{eq:initial_step}
\left(L^{(0)},R^{(0)}\right) = ((3_r,2_r,1_r),(2_\ell,1_\ell,3_\ell)).
\end{equation}
and hence
\begin{equation} \label{eq:initial_step_values}
\begin{array}{l}
\LLL_1^{(0)} = \{(0,0)\} \quad \LLL_2^{(0)} = \{(0,0),(0,0) \tau_h\} \quad \LLL_3^{(0)} = \{(0,0),(0,1)\} \\
\RRR_1^{(0)} = \{(0,0)\} \quad \RRR_2^{(0)} = \{(0,0),(1,0)\} \quad \RRR_3^{(0)} = \{(0,0),(0,0)\tau_v\}.
\end{array}
\end{equation}
The six words $L_i^{(k)}$ and $R_i^{(k)}$ for $i=1,2,3$ are defined recursively by the substitutions $\sigma_\ell$ and $\sigma_r$ in~(\ref{eq:sigma_l_and_sigma_r}). More precisely denoting $L = L^{(k-1)}$, $L' = L^{(k)}$, $R =R^{(k-1)}$ and $R' = R^{(k)}$ we have
\begin{equation} \label{eq:odd_and_even_steps}
\begin{array}{l|l}
\text{for odd steps $k$} & \text{for even steps $k$} \\
\hline
R' = R & L' = L \\
L'_1 = (R_1 R_2)^{m_k}\ L_1 & R'_1 = (L_1 L_3)^{m_k}\ R_1\\
L'_2 = (R_2 R_1)^{m_k}\ L_2 & R'_2 = (L_2)^{n_k}\ R_2\\
L'_3 = (R_3)^{n_k}\ L_3 & R'_3 = (L_3 L_1)^{m_k}\ R_3 
\end{array}
\end{equation}
In order to simplify notations, we use the above notations in many proofs: $R'$ and $R$ (resp. $L'$ and $L$) for $R^{(k+1)}$ and $R^{(k)}$ (resp. for $L^{(k+1)}$ and $L^{(k)}$).

The first step of the proof consists in analyzing the value of the cocycle at the endpoints of each of the pieces $L^{(k)}$ and $R^{(k)}$. In Lemma~\ref{lem:shape}, we prove that the endpoints are always oriented in the same way for all $k$, more precisely the value of the cocycle $g$ with value in $K$ is constant. Then, using this property, we prove Lemma~\ref{lem:growth} which gives an explicit values for theses endpoints.

As it was defined in Section~\ref{subsection:windtree_cocycle}, the wind-tree cocycle $f$ decomposes into two parts. The first one $g$ with values in $K = \ZZ/2 \times \ZZ/2$ and the other one with values in $\ZZ^2$. Let $g: \mathcal{A}^* \rightarrow K$ be the composition of $f$ with the projection $G \rightarrow K$. There is a natural lift of $K$ into $G$ and we set for $w \in \mathcal{A}^*$, $\overline{f}(w) = f(w) g(w) \in \ZZ^2$.
\begin{lemma} \label{lem:shape}
We have for any $k \geq 0$
\[
\begin{array}{lll}
g \left(L_1^{(k)} \right) = g(3_r) = 1 & g \left(L_2^{(k)} \right) = g(2_r) = \tau_h & g \left(L_3^{(k)} \right) = g(1_r) = 1 \\
g \left(R_1^{(k)} \right) = g(2_\ell) = 1 & g \left(R_2^{(k)}\right) = g(1_\ell) = 1 & g\left(R_3^{(k)}\right) = g(3_\ell) = \tau_v
\end{array}.
\]
and
\[
\overline{f} \left(L_1^{(k)} \right) = \overline{f} \left(L_2^{(k)}\right) \in \NN \times \NN
\quad \text{and} \quad
\overline{f}\left(R_1^{(k)}\right) = \overline{f}\left(R_3^{(k)}\right) \in \NN \times \NN
\]
and
\[
\overline{f} \left(L_3^{(k)}\right) \in \{0\} \times \NN
\quad \text{and} \quad
\overline{f} \left(R_2^{(k)}\right) \in \NN \times \{0\}.
\]
\end{lemma}

\begin{proof}
The statement is true for $k=0$ from the definition of $\left(L^{(0)},R^{(0)}\right)$ in~(\ref{eq:initial_step}) and definition of the wind-tree cocycle. Then we proceed by induction. We do the proof for odd steps, the case of even steps being similar. Assume that $k$ is even and that $L = L^{(k-1)}$ and $R = R^{(k-1)}$ satisfies the conclusion of the lemma. Let $L' = L^{(k)}$ and $R'=R^{(k)}$. As $R' = R$ the conclusion holds for $R'$. From the definition of $L'$ we have 
\[
f (L'_1) = f ( ( R_1 R_2)^{m_k}\ L_1 ) = f ( R_1 R_2 )^{m_k}\ f (L_1).
\]
From induction hypothesis, all of $f(R_1)$, $f(R_2)$ and $f(L_1)$ belongs to $\NN^2$. Hence $f(L'_1) \in \NN$. Similarly for $L'_2$ we have
\[
f (L'_2) = f (R_2 R_1)^{m_k}\ f (L_2).
\]
From induction hypothesis, $f(L_2)$ is of the form $(h,v) \tau_h$ with $h,v \in \NN$ and hence get the conclusion for $L'_2$. Now consider the case of $L'_3$. As $n_k = 2 n'$ is even (assumption in theorem~\ref{thm:even_implies_div}) one has
\[
f(L'_3) = (f (R_3)^2 )^{n'}\ f(L_3).
\]
But as $f(R_3)$ is of the form $(h,v) \tau_v$ with $(h,v) \in \NN$, we have $f(R_3)^2 = (0,2v)$. We hence get the conclusion for $R'_3$. This ends the proof of the lemma.
\end{proof}
Now, we compute explicitly the sequence $\overline{f}\left(L_j^{(k)}\right),\overline{f}\left(R_j^{(k)}\right)$ which from Lemma~\ref{lem:shape} depends only on six parameters. Let $X^{(k)} = \left(x_1^{(k)},x_2^{(k)},x_3^{(k)} \right)$ and $Y^{(k)} = \left(y_1^{(k)},y_2^{(k)},y_3^{(k)}\right)$ be the vectors with non negative integer entries such that
\begin{equation} \label{eq:coordinates_x_y}
\begin{array}{l@{\qquad}l}
\overline{f}\left(L_1^{(k)}\right) = \overline{f}\left(L_2^{(k)}\right) = \left(x_1^{(k)},x_2^{(k)}\right) &
\overline{f}\left(L_3^{(k)}\right) = \left(0,x_3^{(k)}\right) \\
\overline{f}\left(R_1^{(k)}\right) = \overline{f}\left(R_3^{(k)}\right) = \left(y_2^{(k)},y_1^{(k)}\right) &
\overline{f}\left(R_2^{(k)}\right) = \left(y_3^{(k)},0\right)
\end{array}
\end{equation}
Our convention for $y_1$ and $y_2$ may seem strange but is explained by the nice formula in the lemma below.
\begin{lemma} \label{lem:growth}
For odd steps, only the coordinates of $X$ are modified as
\[
X^{(2k+1)} = X^{(2k)} + 
\begin{pmatrix}
0 & m_{2k+1} & m_{2k+1} \\
m_{2k+1} & 0 & 0 \\
 n_{2k+1} & 0& 0
\end{pmatrix} Y^{(2k)}.
\]
For even steps, only $Y$ is modified as
\[
Y^{(2k)} = Y^{(2k-1)} +
\begin{pmatrix}
0 & m_{2k} & m_{2k} \\
m_{2k} & 0 & 0 \\
n_{2k} & 0 & 0
\end{pmatrix} X^{(2k-1)}.
\]
\end{lemma}
\begin{proof}
We omit the proof which proceeds by induction and follows the one of Lemma~\ref{lem:shape}.
\end{proof}
We now build explicit ``boxes'' around the trajectory. More precisely we find the minimum and maximum values of each coordinates of the sets $\LLL_i^{(k)}$ and $\RRR_i^{(k)}$. In order to take care of the horizontal excursions of $L_3^{(k)}$ and vertical excursions of $R_2^{(k)}$, we add one coordinate to the vectors $X'$ and $Y'$. Let $x_4^{(0)} = y_4^{(0)} = 1$ and define recursively for odd steps
\begin{equation} \label{eq:fourth_coordinate_odd}
x_4^{(2k+1)} = \left\{
\begin{array}{ll}
\max \left(x_4^{(2k)}, y_2^{(2k)} \right) & \text{if $n_{2k+1} \not= 0$,} \\
x_4^{(2k)} & \text{otherwise},
\end{array} \right.
\qquad \text{and} \qquad
y_4^{(2k+1)} = y_4^{(2k)}
\end{equation}
and for even step
\begin{equation} \label{eq:fourth_coordinate_even}
x_4^{(2k)} = x_4^{(2k-1)} \qquad \text{and} \qquad
y_4^{(2k)} = 
\left\{ \begin{array}{ll}
\max \left(y_4^{(2k-1)}, x_2^{(2k-1)} \right) & \text{if $n_{2k} \not= 0$,} \\
y_4^{(2k-1)} & \text{otherwise}
\end{array} \right.
\end{equation}
We first start by the formal definition of a box.
\begin{definition}
Let $\pi_h:\ZZ^2 \rightarrow \ZZ$ (resp. $\pi_v: \ZZ^2 \rightarrow \ZZ$) be the projection on the first (resp. the second) coordinate. Let $A \subset \ZZ^2$. The \emph{box} of $A$ is the $4$-tuple 
\[
\Box(A) = (\min \pi_h(A), \min \pi_v(A), \max \pi_h(A), \max \pi_v(A)) \in \ZZ^4.
\]
\end{definition}
By extension, we call for $i=1,2,3$ the \emph{box} of the word $L_i^{(k)}$ (resp. $R_i^{(k)}$) the box of the subset $\LLL_i^{(k)}$ (resp. $\RRR_i^{(k)}$). The boxes around the pieces $(L^{(k)},R^{(k)})$ is given in the following lemma.
\begin{lemma} \label{lem:boxes}
We have
\[
\begin{array}{l@{\qquad}l}
\multicolumn{2}{c}{\Box \left(L_1^{(k)}\right) = \Box \left(L_2^{(k)}\right) = \left(0,0,x_1^{(k)},x_2^{(k)}\right)},
\\
\multicolumn{2}{c}{\Box \left(R_1^{(k)}\right) = \Box \left(R_3^{(k)}\right) = \left(0,0,y_2^{(k)},y_1^{(k)}\right)},
\\
\Box \left(L_3^{(2i)}\right) = \left(0,0,x_4^{(k)},x_3^{(k)}\right),
&
\Box \left(R_2^{(2i)}\right) = \left(0,0,y_3^{(k)},y_4^{(k)}\right).
\end{array}
\]
\end{lemma}

\begin{proof}
The cases of $L_1^{(k)}$, $L_2^{(k)}$ is straightforward from the proof of Lemma~\ref{lem:shape} as well as the case of $R_1^{(k)}$ and $R_3^{(k)}$. We prove by induction the formula for the box of $L_3^{(k)}$, the case of $L_3^{(k)}$ being similar. There is nothing to prove at even steps as $L_3^{(2k)} = L_3^{(2k-1)}$. Assume that the conclusion of the lemma holds at an even step $k-1$ and denote $L=L^{(k-1)}$ and $L'=L^{(k)}$. We recall that
\[
L'_3 = (R_3 )^{n_k}\ L_3
\]
where $n_k$ is an even number by assumption in Theorem~\ref{thm:even_implies_div}. If $n_k = 0$, then $L'_3 = L_3$, and from~(\ref{lem:growth}) we have $x'_3 = x_3$ and from~(\ref{eq:fourth_coordinate_odd}) $x'_4 = x_4$. Hence the box fits in this case. Now assume that $n_k \not= 0$. We know from our induction hypothesis that $\Box(R_3) = (0,0,y_2,y_1)$. As the word $R_3$ ends with $\tau_v$, for any $n \geq 1$ we have $\Box\left((R_3)^n\right) = \left(0,0,y_2, n\, y_1\right)$. From Lemma~\ref{lem:growth} we have $f\left(R'_3\right)^{(n_k)} = (0,n_k\, y_1) \in \{0\} \times \NN$ and $f\left(L'_3\right) = (0,x_3)$. The word $L'_3$ is the concatenation of $(R_3)^{n_k}$ and $L_3$ and is hence contained in the box with bottom-left corner $(0,0)$ and up-right corner $\left(\max \left(x_4, y_2 \right), n_k\, y_1 + x_3\right) = \left(x'_4 , x'_3\right)$.
\end{proof}

The following lemma states that the trajectories in the billiard are self-avoiding at large scales. In other words if the trajectory crosses at time $t_0$ and $t_1$ then the difference $|t_1-t_0|$ should be small. In the proof of Theorem~\ref{thm:even_implies_div} below, we refine the argument to prove that at small scales intersection does not appear as well.
\begin{lemma} \label{lem:self_avoiding}
Let $k \geq 2$ be such that all entries of $X^{(k-1)}$ and $Y^{(k-1)}$ are positive. Let $W'$ be one of the six words $L_i^{(k)}$ or $R_i^{(k)}$ for $i=1,2,3$ and $W' = W_1 W_2 \ldots W_p$ its decomposition given by the Ferenczi-Zamboni induction where each $W_j$ equals one of the six words $L_i^{(k-1)}$ or $R_i^{(k-1)}$ for $i=1,2,3$. We denote by $\mathcal{W}_j$ the set of values taken by the wind-tree cocycle on $W_j$. Let $j$ and $j'$ be two distinct elements of $\{1,\ldots,p\}$. Then the subsets $f(W_1 \ldots W_{j-1})\ \mathcal{W}_j$ and  $f(W_1 \ldots W_{{j'}-1})\ \mathcal{W}_{j'}$ are disjoint if $|j-j'| \not= 1$ and have only one intersection point otherwise.
\end{lemma}

\begin{proof}
The decomposition of $W'$ depends on the parity of $k$ and is given by the rules~(\ref{eq:odd_and_even_steps}). We do the proof at an odd step of the induction. Let $W'= L'_1 = (R_1 R_2)^{m_k} L_1 = W_1 W_2 \ldots W_{2m_k+1}$ where as in the proof of the preceding lemmas $L'_i = L^{(k)}_i$, $R'_i = R^{(k)}_i$, $L_i = L^{(k-1)}_i$ and $R_i=R^{(k-1)}_i$. From lemmas~\ref{lem:growth} and~\ref{lem:boxes}, we know that the values of the wind-tree cocycle on $R_1$ and $R_2$ respectively ends in the top-right corner and the bottom-right corner of the respective box $\Box(R_1)$ and $\Box(R_2)$. From Lemma~\ref{lem:shape}, we deduce that in the word $W' = (R_1 R_2)^{m_k} L_1$ each of the individual box $f(W_1 \ldots W_{j-1})\ \mathcal{W}_j$ is glued to the preceding at only one point which is $f(W_1 \ldots W_{j-1})$. The proof works the same for $L'_2 = (R_2 R_1)^{m_k} L_2$ and $L'_3 = (R_3)^{n_k} L_3$.

Because of our assumption on the entries of $X^{(k-1)}$ and $Y^{(k-1)}$, the intersection of two non consecutive boxes is empty.
\end{proof}

We are now ready to prove our main theorem.
\begin{proof}[Proof of Theorem~\ref{thm:even_implies_div}]
Let $w \in \mathcal{A}^\ZZ$ be a cutting-sequence of an orbit of the interval exchange determined by $(a,b,\theta)$ which satisfies the hypothesis of Theorem~\ref{thm:even_implies_div}. Then, for each $k$, we decompose $w$ as a word on the alphabet $\mathcal{A}^{(k)} = \{L^{(k)}_1, L^{(k)}_2,L^{(k)}_3,R^{(k)}_1,R^{(k)}_2,R^{(k)}_3\}$. We stress that the origin of $w$ (the letter in position $0$) should be shifted in order to have a decomposition on $\mathcal{A}^{(k)}$ which starts at position $0$. But this does not matter for our purpose and consider that $w$ has no origin when it is decomposed with respect to the alphabet $\mathcal{A}^{(k)}$. From Lemma~\ref{lem:self_avoiding}, we know that for $k$ big enough, the pieces of size $k$ contained in a piece of size $k+1$ are disjoint except for one possible value which occurs at the end of a box and the begining of the next one. Hence, if the box that contains the origin grows arbitrarily on the left and on the right, the trajectory determined by $w$ is divergent. But if the box does not growth arbitrarily on the right, say it is blocked at $N$, then it means that in the future the orbit of $p$ encounters a singularity of the interval exchange transformation after $N$ steps.

We now refine the argument at small scales to prove that the trajectory in the infinite billiard is self-avoiding. It follows from the combinatorics of the surface $\Lsurf(a,b)$, that two consecutive pieces of level $k$ may intersect in only few cases (see Figure~\ref{fig:Lab_geodesic}) which corresponds to block on which the cocycle remains constant:
\begin{itemize}
\item either in the word $w = 3_r 3_\ell^n$ where $n \geq 0$,
\item or in the form $w = 2_\ell 2_r^n$ where $n \geq 0$.
\end{itemize}
In the first case (resp. the second) the word is followed by $1_r$ (resp. $1_\ell$). In the billiard table, the cutting sequence $w$ lift to a piece of trajectory which reflects between two horizontally (resp. vertically) consecutive scatterers but does not reflect vertically (resp. horizontally). In particular, the trajectory is self-avoiding in $w$.
\end{proof}

\section{An example}
\begin{figure}
\begin{center}
\includegraphics[scale=0.37]{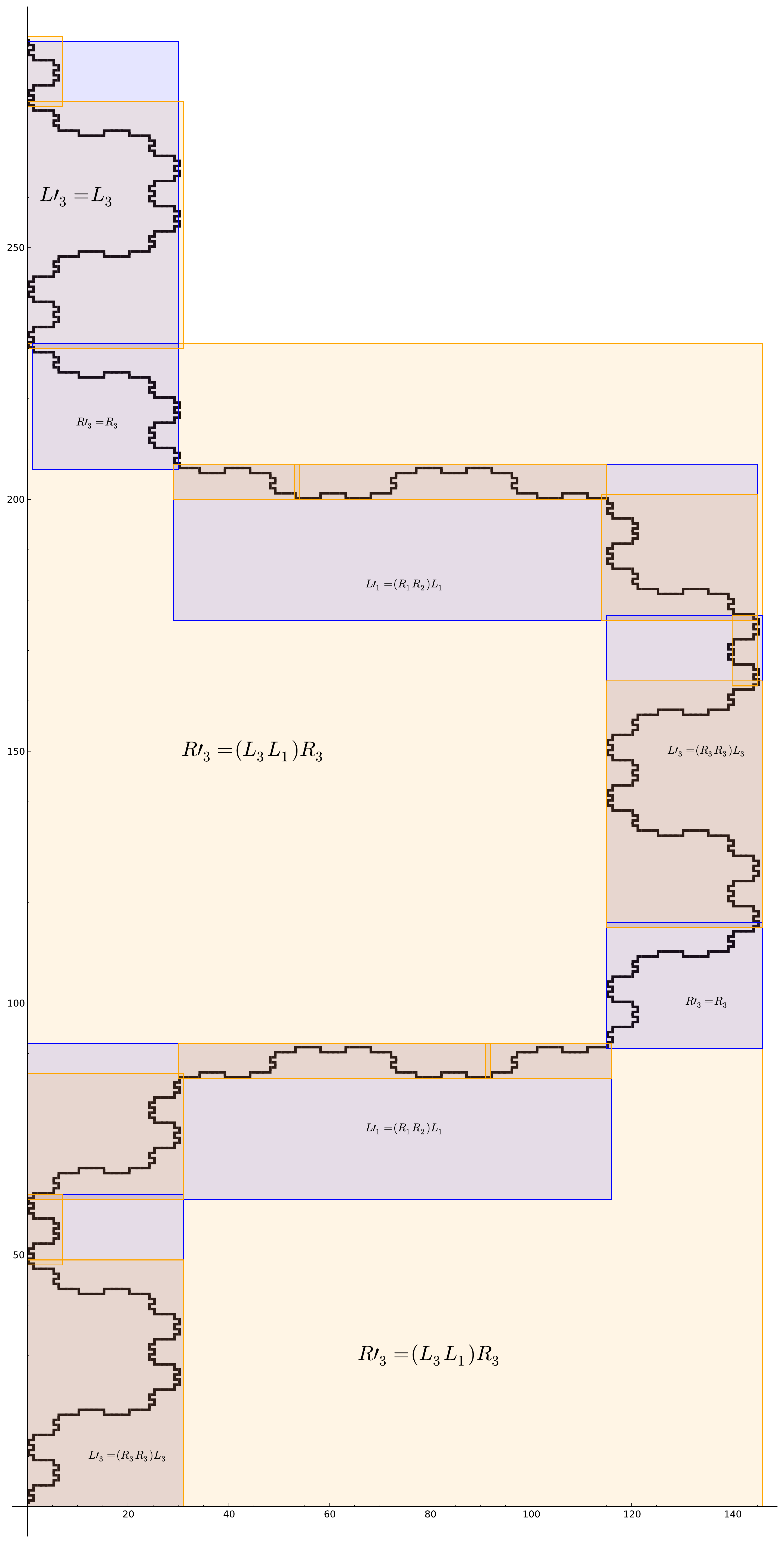}
\end{center}
\caption{Three levels of boxes for the wind-tree cocycle over $\Lsurf(1/2,1/2)$ and the slope $\sqrt{2}-1$. The line in black corresponds to the quadrilaterals visited by an orbit in the wind-tree model $T(1/2,1/2)$.}
\label{fig:L22_sqrt2p1}
\end{figure}

We now consider the example of the periodic sequence of convergents $((1,2),(1,2),\ldots)$ associated to the square tiled surface $\Lsurf(1/2,1/2)$ and the slope $\sqrt{2}-1$. For this slope, the Ferenczi-Zamboni induction is periodic or in other words the interval exchange transformation is self-similar. The Figure~\ref{fig:L22_sqrt2p1} shows a three level of boxes for an associated orbit in the infinite billiard table $T(1/2,1/2)$.

\newpage

\end{document}